\documentclass[12pt,reqno]{amsart}
\usepackage{type1ec}
\usepackage[utf8]{inputenc}
\usepackage[T1]{fontenc}
\usepackage{amssymb,amsfonts,amsthm}
\usepackage{graphicx}
\usepackage{tikz}              
\usetikzlibrary{calc}  
\usepackage{amsthm}
\usepackage{amssymb}
\usepackage[francais,english]{babel}
\usepackage{setspace}
\usepackage{hyperref}

\usepackage{color}

\definecolor{Red}{cmyk}{0,1,1,0.2}

\usepackage{amsfonts,amsmath,layout}

\newcommand{\Q}{\mathbb Q}

\newcommand{\R}{\mathbb R}

\def\R{\mathbb R}

\def\E{\mathbb E}
\def\P{\mathbb P}

\newcommand{\be}{\begin{equation}}
\newcommand{\ee}{\end{equation}}

\def\1{{\bf 1}}

\def\ds{\displaystyle}

\newtheorem{Theorem}{Theorem}[section]
\newtheorem{Definition}[Theorem]{Definition}
\newtheorem{Proposition}[Theorem]{Proposition}
\newtheorem{Lemma}[Theorem]{Lemma}
\newtheorem{Corollary}[Theorem]{Corollary}
\newtheorem{Remark}[Theorem]{Remark}

\newtheorem{Sketch of proof}[Theorem]{Sketch of proof}
\begin{document}
\title[Stochastic optimal scattering control]{Stochastic scattering control of spider diffusion governed by an optimal diffraction probability measure selected from its own local time}
\author[Isaac Ohavi]{Isaac Ohavi$^{\aleph}$\\
\\{$^\aleph$Hebrew University of
Jerusalem, Department of Mathematics and Statistics, Israël}}
\email{isaac.ohavi@mail.huji.ac.il \& isaac.ohavi@gmail.com}
\thanks{This research  was supported by the GIF grant 1489-304.6/2019.\\
I am grateful to Miguel Martinez$^{\beth}$ for the several discussions and collaborations on spider diffusions and their related partial differential equations. I am also grateful to Ioannis Karatzas $^\gimel$, and the different members of the seminary: Probability and Control theory - Columbia University, New York.\\
$^\beth$ Université Gustave Eiffel, LAMA UMR 8050, France. E-mail: miguel.martinez@univ-eiffel.fr\\
$^\gimel$ Department of Mathematics, Columbia University, New York, NY 10027, USA. E-mail: ik1@columbia.edu.}
\dedicatory{Version: \today}
\begin{abstract}The purpose of this article is to study a new problem of stochastic control, related to Walsh's spider diffusion, named: {\it stochastic optimal scattering control}. The optimal scattering control of the spider diffusion at the junction point is governed by an 
appropriate and highly non-trivial condition of the “Kirchhoff Law” type, involving an optimal diffraction probability measure selected from the own local time of the spider process at the vertex. In this work, we prove first the weak dynamic programming principle in the spirit of \cite{ElKaroui}, adapted to the new class of spider diffusion introduced recently in \cite{Spider}-\cite{Spider 2}. Thereafter, we show that the value function of the problem is characterized uniquely in terms of a Hamilton Jacobi Bellman (HJB) system posed on a star-shaped network, having a new boundary condition at the vertex called : {\it non linear local-time Kirchhoff's transmission}. The key main point is to use the recent comparison theorem obtained in \cite{Ohavi Walsh PDE} that has significantly unlocked the study of this type of problem.  We conclude by discussing the formulation of stochastic scattering control problems, where there is no dependency w.r.t. the local-time variable, for which their well-posedness appear as a simpler consequence of the results of this work and the advances contained in \cite{Ohavi Walsh PDE}. 
\end{abstract}
\maketitle
{\small \textbf{Key words:}  Stochastic scattering control, DPP, HJB equations and non linear local-time Kirchhoff's boundary condition, comparison principle and viscosity solutions, local-time.\\
\textbf{MSC2020 subject classifications:}  34B45, 35R02, 49L20, 49L25, 60G46, 78A45, 93E20.}
\section{Introduction}\label{sec:intro}
\subsection{General Introduction:}
Walsh's spider process whose spinning measure and coefficients are allowed to depend on the local time at the junction vertex and the running time, were introduced by Martinez-Ohavi in \cite{Spider}. In this contribution, the authors construct the process and prove that there is uniqueness in the weak sense of its law on the joint spider's path space-local time space. 

Up to our knowledge, the last contribution \cite{Spider} states the first result in literature that provides existence and uniqueness of a Walsh spider's diffusion with non constant spinning measure. More precisely, let $I$ a positive integer ($I\geq 2$). Denote $\mathcal{N}$ the following star-shaped network:
\begin{align*}
\mathcal{N}~~:=~~\{{\bf 0}\}\cup \big((0,\infty)\times [I]\big),~~[I]:=\{1 \ldots I\}.
\end{align*} 
where ${\bf 0} = \{(0,j), j\in [I]\}$ is the junction vertex equivalence class. Assume we are given $I$ pairs $(\sigma_{i},b_{i})_{i\in [I]}$ of mild coefficients from $[0,T]\times[0,+\infty)^2$ to $\R$ satisfying the following condition of ellipticity: $\forall i\in [I],~\sigma_i>0$. We are also given $\big(\mathcal{S}_1,\ldots,\mathcal{S}_I)\in \mathcal{C}\big([0,T]\times[0,+\infty),\R_+\big)^I$ a positive spinning measure (depending on the running time and the local-time) such that:
$$\forall (t,l)\in [0,T]\times[0,+\infty),~~\displaystyle \sum_{i=1}^I \mathcal{S}_i(t,l)=1.$$ 
The statement of the martingale problem $\big(\mathcal{S}_{pi}-\mathcal{M}_{ar}\big)$ studied in \cite{Spider}, is formulated for continuous process $(x(s),i(s),l(s))_{s\in [0,T]}$ with state space $\mathcal{N}\times \R^+$, starting at a point $(x_\ast,i_\ast,l_\ast)\in \mathcal{N}\times [0,+\infty)$ at time $t$, such that for any $f\in\mathcal{C}^{1,2,1}_b\big([0,T]\times \mathcal{N}\times [0,+\infty),\R\big)$:
\begin{align} \label{eq: def Martingale Problem}
&\nonumber \Bigg(f_{i(s)}(s,x(s),l(s))- f_{i_\ast}(t,x_\ast,l_\ast)\displaystyle\\
&\nonumber - \int_{t}^{s}\Big(\partial_tf_{i(u)}(u,x(u),l(u))+\displaystyle\frac{1}{2}\sigma_{i(u)}^2(u,x(u),l(u))\partial_{xx}^2f_{i(u)}(u,x(u),l(u))\\
&\nonumber +b_{i(u)}(u,x(u),l(u))\partial_xf_{i(u)}(u,x(u)l(u))\Big)du\\
&-\int_{t}^{s}\Big(\partial_lf_{i(u)}(u,0,l(u))+\displaystyle\sum_{j=1}^{I}\alpha_j(u,l(u))\partial_{x}f_{j}(u,0,l(u))\Big)dl(u)~\Bigg)_{t\leq s\leq T},
\end{align} 
is a martingale under some probability measure $\P_t^{(x_\ast,i_\ast,l_\ast)}$. It was deduced also that $(\ell^0_s:= l(s)-l_\ast)_{s\in [t,]}$ is in fact the local time of the spider $(x(s),i(s))_{s\in [t,]}$ at $\bf 0$, increasing only on the times when the component $(x(s))_{s\in [t,]}$ equals $0$. 

The existence proof relies on taking as a starting point first the result dealing with the weak existence of spider diffusion with constant spinning measure, introduced for the first time in \cite{freidlinS}-\cite{Freidlin-Wentzell-2}. Thereafter, the authors in \cite{Spider} constructed 'by hand' a solution to the last martingale problem in order to take the presence of the local time in all the leading coefficients into account. Of course, since the local time was added in the picture, the canonical space had to be extended accordingly. More precisely, the construction was performed with a careful adaptation of the seminal construction for solutions of classical martingale problems that have $\R^d$ as the underlying state space and combined concatenation of probability measures with tension arguments, such as introduced in the essential reference book \cite{Stroock}. 

The uniqueness proof is more involved and comes essentially from the results obtained  in \cite{Martinez-Ohavi EDP}, dealing with the well-posedness of the associated 'new' parabolic operator, generating the martingale problem $\big(\mathcal{S}_{pi}-\mathcal{M}_{ar}\big)$ \eqref{eq: def Martingale Problem}.
More precisely the following (forward) system posed on the domain $[0,T]\times \mathcal{N}\times [0+\infty)$ :
\begin{eqnarray}\label{eq : pde with l}
\begin{cases}
\textbf{Linear parabolic equation parameterized}\\
\textbf{by the local-time on each ray:}\\
\partial_tu_i(t,x,l)-a_i(t,x,l)\partial_x^2u_i(t,x,l)
+b_i(t,x,l)\partial_xu_i(t,x,l)\\
\hspace{0,4 cm}+c_i(t,x,l)u_i(t,x,l)=f_i(t,x,l),~~(t,x,l)\in (0,T)\times (0,+\infty)^2,\\
\textbf{Linear local-time Kirchhoff's boundary condition at }\bf 0:\\
\partial_lu(t,0,l)+\displaystyle \sum_{i=1}^I \mathcal{S}_i(t,l)\partial_xu_i(t,0,l)=\phi(t,l),~~(t,l)\in(0,T)\times(0,+\infty),\\
\textbf{Continuity condition at } \bf 0:\\
\forall (i,j)\in[I]^2,~~u_i(t,0,l)=u_j(t,0,l)=u(t,0,l),~~(t,l)\in|0,T]\times[0,K],\\
\textbf{Initial condition:}\\
\forall i\in[I],~~ u_i(0,x,l)=g_i(x,l),~~(x,l)\in[0,+\infty)^2,
\end{cases}
\end{eqnarray}
has been studied recently by Martinez-Ohavi in \cite{Martinez-Ohavi EDP}. From a PDE technical aspect, since the variable $l$ drives dynamically the system only at the junction point $\bf 0$ with the presence of the derivative $\partial_lu(t,0,l)$ in the local time Kirchhoff's boundary condition, the main challenge was to understand the regularity of the solution. Under mild assumptions, it is shown  that classical solutions of the system \eqref{eq : pde with l} belong to the class $\mathcal{C}^{1,2}$ in the interior of each edge and $\mathcal{C}^{0,1}$ in the whole domain (with respect to the time-space variables $(t,x)$). One can expect a regularity in the class $\mathcal{C}^1$ for $l\mapsto u(t,0,l)$ and this is indeed the case (see Theorem 2.4 and point $iv)$ in Definition 2.1 in \cite{Martinez-Ohavi EDP}). Another technical aspect, was to obtain an H\"{o}lder continuity of the partial functions $l\mapsto \Big(\partial_tu_i(t,x,l),\partial_xu_i(t,x,l),\partial_x^2u_i(t,x,l)\Big)$ for any $x>0$. It is shown, using Schauder's estimates, that such regularity is guaranteed by the central assumption on the ellipticity of the diffusion coefficients on each ray, together with the mild dependency of the coefficients and free term with respect to the variable $l$. Finally, note that it has also been proved that the derivative with respect to the variable $l$: $\partial_lu_i(t,x,l)$ was locally uniformly bounded in the whole domain. The comparison theorem holds true under such regularity class; see Theorem 2.6 in \cite{Martinez-Ohavi EDP}. Remark that the results contained in \cite{Martinez-Ohavi EDP} extend those obtained by Von Below in \cite{Von-Below}, which were -- up to our knowledge -- the only reference on the well-posedness of heat equations on graphs with time variable coefficients and classical Kirchhoff's condition:
\begin{eqnarray*}
\sum_{i=1}^I \mathcal{S}_i(t)\partial_x f_i(t,0)=0.
\end{eqnarray*}
Regarding the last arguments, we see that the class of regularity of the solution of system \eqref{eq : pde with l} is then weaker than the class $\mathcal{C}^{1,2,1}_b([0,T]\times \mathcal{N}\times [0,+\infty),\R)$, used for the formulation of the spider martingale problem $\big(\mathcal{S}_{pi}-\mathcal{M}_{ar}\big)$ \eqref{eq: def Martingale Problem}. If one wants to achieve the proof of the weak uniqueness of $\big(\mathcal{S}_{pi}-\mathcal{M}_{ar}\big)$ (Theorem 3.1 in \cite{Spider}), it appears that we have to show that the martingale formulation also holds true under the class of regularity of the solution of system \eqref{eq : pde with l}. This was achieved in \cite{Spider} Proposition 6.3, with the aid of a  special regularization procedure that is not standard, which ensures continuity at the junction vertex $\bf 0$.

The companion article of \cite{Spider}, also written by the same authors in \cite{Spider 2}, derives some crucial properties of this type of new spider process. Therein, among them: the Markov's property, It\^o's formula, some Feynmann Kac's representations for backward systems of type \eqref{eq : pde with l}, the behavior of the density of the process $x$ at $\bf 0$ and some approximations of the local-time. Moreover, \cite{Spider 2} gives an interpretation of the probability coefficients of diffraction $(\mathcal{S}_i)_{i \in [I]}$. Assume that $(x,i)=\bf 0$. Let $\P_t^{\bf 0,\ell}$ be the solution of the $\big(\mathcal{S}_{pi}-\mathcal{M}_{ar}\big)$ martingale problem, starting at time $t$ from the junction point $(x,i)=\bf 0$, with a local-time level equal to $\ell$. Then, for any $\delta>0$, if we introduce the following stopping time:
$$\theta^\delta:=\inf\big\{~s\ge 0,~~x(s)=\delta~\big\},$$
we have:
\begin{eqnarray}\label{eq distr instant diffra temp local}
\forall i\in [I],~~\lim_{\delta \searrow 0} \P_t^{\bf 0,\ell}\big(~i(\theta^\delta)=i~\big) =\mathcal{S}_i(t,\ell).  
\end{eqnarray} 
This result allows us to state that as soon as the spider process $(x,i)$ reaches the junction point $\bf 0$ at time $t$, with a level of local time $\ell$, then the 'instantaneous' probability distribution for $(x,i)$ to be scattered along the ray $\mathcal{R}_i$ is then equal to the corresponding spinning coefficient $\mathcal{S}_i(t,\ell)$. 

The aim of this work is: \textbf{Formulate a problem of scattering control in the finite horizon time framework, and characterize its optimal diffraction probability  $\mathcal{S}^*_i(t,\ell)$ selected from the own local-time at $\bf 0$, with the aid of the corresponding value function that is the unique solution of a HJB system posed on $\mathcal{N}$}.

As soon as one wants to characterize uniquely the value function of a stochastic control problem, it is well known that non linear PDE and in particular viscosity solution theory plays a determinant role. Let us describe the last advances obtained in \cite{Ohavi Walsh PDE}, that are, to our knowledge, the first reference in viscosity theory, appropriated to a stochastic scattering problem for a spider, with optimal probability of diffraction selected from its own local-time.

To prove the comparison theorem for a second order HJB system posed on star-shaped network having the non linear local-time Kirchhoff's condition, the author built in \cite{Ohavi Walsh PDE} test functions at the junction point $\bf 0$ solutions of ODE with coefficients that may be viewed as a kind of envelope of all possible errors of the {\it speed of the Hamiltonians}. The key point in the construction was to impose a local-time derivative with respect to variable $l$, at $\bf 0$ - $\partial_l\phi(0,l)$ - that absorb all the error term induced by - the {\it Kirchhoff's speed of the Hamiltonians}. This was the first result of uniqueness for HJB elliptic PDE system posed on a star-shaped-network, having a nonlinear Kirchhoff's condition and non vanishing viscosity at the vertex. Let us emphasize that in \cite{Ohavi Walsh PDE} the introduction of the external deterministic ’local-time’ variable $l$ -- that is, the counterpart to the local time $\ell$ -- is one of the crucial ingredients to obtain the comparison principle. Note that even without the presence of the external variable $l$ in the HJB problems already studied in the literature, the 'artificial' introduction of this external variable in the problem allowed to extend the main results contained in \cite{Lions Souganidis 1} and \cite{Lions Souganidis 2} to the fully non linear and non degenerate framework. On the other hand, it is important also to emphasize that the comparison theorem for viscosity solutions in \cite{Ohavi Walsh PDE} is stated in the strong sense for Kirchhoff's boundary condition at $\bf 0$, namely without any dependency of the values of the Hamiltonians at $\bf 0$, which is also innovative for Neumann problems in the non linear case.

To the best of our expertise, the characterization of the value function of a problem of control involving a spider diffusion, could be achieved in the past literature, only via a verification theorem. The well-posedness of quasi linear systems (elliptic and parabolic) was studied in \cite{Ohavi PDE} and thereafter used recently, for example, in MFG problems posed on networks in \cite{Berry} for sticky spider diffusion.

Let us now make an important remark. The difficulty when studying spider processes comes from the fact that the natural filtration generated by a spider process is not a Brownian filtration as soon as the underlying spider web possesses three branches or more. As mentioned in \cite{Lejay}: {\it ''this fact has given rise to an abundant literature on Brownian filtrations''}. We mention \cite{DeMeyer} for a comprehensive argument that the natural filtration cannot be generated by a Brownian motion when $I\geq 3$ and also the unavoidable reference \cite{Barlow} (see also \cite{tsirel}) for a theoretical study of these questions.\\
Hence, as soon as one wants to formulate a stochastic control problem for spider processes, the dynamic programming principle can be established, therefore only in the weak sense, as soon as $I\ge 3$. If no controls appear at the junction point $\bf 0$ in the diffraction terms governing the Kirchhoff's condition, we believe that the arguments introduced in \cite{ElKaroui} can be easily adapted. But as soon the controls appear at $\bf 0$, one has to carefully be able to adapt the results of \cite{ElKaroui}, and has to prove especially the compactness of the admissible rules at $\bf 0$. We will see that another fundamental result of this contribution will be to properly prove the dynamic programming principle, as soon as control terms appear at the vertex.
This contribution being already quite long, note that we will only treat the strictly elliptical case, in order to alleviate in particular the characterization of the function. The dynamic programming principle could be extended to to degenerate case, as it has been already done in \cite{ElKaroui}.

We believe that some expected applications, for example in physics could be undertaken. Once again, so as not to make the reading of this work cumbersome, we describe just some lines of this potential applications in this Introduction, since other contributions could be dedicated to describe concrete examples, and also treat others problems related to stochastic scattering optimal control theory, like stopping time/infinite horizon time problems, existence of optimal feedback control... (that we do not discuss here). The study of optimal diffraction control problems for stochastic Walsh diffusions is 
expected to have implications when one tries to study the diffusive behavior of particles subjected to scattering (or diffraction), for which little physical understanding exists
at present. The theory of quantum trajectories states that quantum systems can be modeled as scattering processes ; we refer the reader to \cite{Scatterin Theory} for more appropriate
details. Problems of optimal light scattering have gained more and more interest for their importance in advanced photonic technologies such as, for example, on-chip interconnects, bio-imaging, solar-cells, or heat-assisted magnetic recording.
In order to illustrate the issues addressed in this paper, let us imagine that a punctual source of light crosses a plane at some point $O$. We constrain a Brownian particle to
move along a finite number of rays with different magnetic and electronic properties, that are joined on a ’spider web’ whose central vertex lies at $O$. When passing at the
vertex junction, the Brownian particle gets directly hit by the punctual source light, and this modifies its electronic properties. Thus, the particle is instantly attracted
in a more privileged manner towards some particular rays of the
spider web and these change according to its modified electron affinity received instantly from the punctual
light source. Since the particle gets directly affected by the time it has spent under the light i.e. the ’local time’ spent by the Brownian particle under the light source of the
vertex, this ’local time’ has a direct influence on the privilege instantaneous directions elected by the particle (i.e. the spinning measure of the motion). In turn, such a device would give information on the light scattering of a punctual light source by a single particle. Note that such a device could not be set up by using a classical planar Brownian
motion particle : because the trajectories do not have bounded variations, it does not seem feasible to pursue the particle with a point laser, and such a motion would never return exactly under a fixed punctual light source point.
The literature in Physics on the subject of scattering light from Brownian particles is extended and technical and we do not pretend at all to know it, even less to
be specialists : this is why we do not prefer to venture into giving a representative and relevant sample here. However, the reader will surely be interested to find in \cite{Clark} what seems to us to be the origins of the study of the relationships between light scattering and Brownian particle motions.

\subsection{Review of literature:} To finish this Introduction, let us give an account of the main works that have been done in literature on spider diffusions, control theory on stratified domains, MFG systems posed on networks, and their related PDE.

The idea was first introduced by Walsh in the epilogue of \cite{Walsh}. As explained in \cite{Barlow-Pitman-Yor}: 
\begin{quote}
Started at a point $z$ in the plane away from the origin $0$, such a process moves like a one-dimensional Brownian motion along the ray joining $z$ and $0$ until it reaches $0$. Then the process is kicked away from $0$ by an entrance law that makes the radial part of the diffusion a reflecting Brownian motion, while randomizing the angular part.
\end{quote} The difficulty arises since, as Walsh explains with certain sense of humor (quoted from \cite{Barlow-Pitman-Yor} and \cite{Walsh}):
\begin{quote}
(...) It is a diffusion that, when away from the origin, is a Brownian motion along a ray, but which has what might be called a round-house singularity at the origin: when the process enters it, it, like Stephen Leacock’s hero, immediately rides off in all directions at once.
\end{quote}
Rigorous construction of Walsh's Brownian motion may be found in \cite{Barlow-Pitman-Yor} (see also the references given therein).

Generalizing the idea of Walsh, diffusions on graphs were introduced in the seminal works of Freidlin and Wentzell \cite{Freidlin-Wentzell-2} and Freidlin and Sheu \cite{freidlinS} for a star-shaped network (and afterwards for open books in \cite{Freidlin}). 
For the sake of conciseness and also because of the difficulty of the subject, in this contribution we will focus on the case where the diffusion lives on a star-shaped network:~these diffusions are called {\it spiders} and, sticking with this metaphor, the state-space junction network may be referred to as {\it the spider's web}. 

Walsh spider diffusion processes are currently being thoroughly studied and extended to various settings. Let us mention the following recent articles among the vast literature on the subject: in \cite{Ichiba} the authors propose the construction of stochastic integral equations related to Walsh semi martingales, in \cite{Ichiba-2} the authors compute the possible stationary distributions, in \cite{Kara control} the authors investigate stopping control problems involving Walsh semi martingales, in \cite{Atar} the authors study related queuing networks, whereas \cite{Bayraktar} addresses the problem of finding related stopping distributions. We refer to the introduction of \cite{Ichiba} for a comprehensive survey, the reader may also find therein many older references on the subject.

Recall that the difficulty when constructing spider processes comes from the fact that the natural filtration generated by a spider process is not a Brownian filtration as soon as the underlying spider web possesses three branches or more. Although the motion behaves as a one dimensional semimartingale during its stay on a particular branch and is driven by a Brownian motion along the branch, the fact that the process {\it 'has more than two directions to choose infinitesimally when moving apart from the junction point'} makes it impossible for a spider diffusion to be adapted to a Brownian filtration (whatever the dimension of the underlying Brownian motion). This crucial fact imposes that it is not possible to apply directly It\^o's stochastic calculus theory to both components $(x,i)$ of the spider process (where $x$ stands for the distance to the vertex junction and $i$ stands for the label of the branch), but only to the distance component $x$.

The relationship between Walsh diffusion processes and skew diffusions (that correspond to the particular case where $I=2$) may be found in the survey \cite{Lejay}. It is notable that skew diffusions solve stochastic differential equations that involve the local time of the unknown process (in this case the natural filtration of the process is a Brownian filtration). 

Although difficult, several constructions of Walsh's diffusions have been proposed in the literature, see  for e.g. \cite{Barlow-Pitman-Yor} for a construction based on Feller’s semigroup theory, \cite{Salisbury} for a construction using the excursion theory for right processes, and also the very recent preprint \cite{Bayraktar-2} that proposes a new construction of Walsh diffusions using time changes of multi-parameter processes. Note that in all these constructions, the spinning measure of the process -- that is strongly related somehow to {\it 'the way of selecting infinitesimally the different branches from the junction vertex'} -- remains constant through time. Hence, to formulate a problem of scattering control for a spider, governed by an optimal diffraction probability measure at the junction point $\bf 0$, it is fundamental to be able to build optimal feedback in the coefficient of diffraction. This was the aim of the contribution \cite{Spider}, that has given the first result in literature dealing with the existence and uniqueness of Walsh processes, with a random non constant spinning measure. Therein, more precisely, the spinning measure depends on the own local time of the process spent at the junction, together with the current running time. 

This contribution deals, to the best of our expertise with the first result in stochastic control theory for reflected processes, having optimal controls at the boundary depending on the local-time. 

For some formulations of stochastic control problems involving Walsh's spider diffusion, we can refer to the recent works in \cite{Kara control}, where an optimal stochastic stopping control problem for a Walsh's planar semi martingale has been studied, but without optimal scattering control at the junction point, since it is assumed that the process is “immediately dispatched along some ray” when it reaches the origin. On the other hand, references on MFG systems posed on networks can be found in \cite{Acdou 1}, \cite{Acdou 2} or \cite{Berry} and the references therein. Note also that boundary stochastic control problems have been studied for classical reflected diffusion in \cite{Bouchard}.

Regarding HJB equations posed on networks, we refer to the recent monograph \cite{Barles book}, that presents the most recent developments in the study of Hamilton-Jacobi Equations and control problems with discontinuities (see Part III for the case of problems on networks). Note that recently, non local first order HJB equations posed on networks, have been studied in \cite{Barles}. 
In \cite{Lions Souganidis 1}, the authors introduce a notion of state-constraint viscosity solutions for one dimensional “junction” - type problems for first order Hamilton-Jacobi equations with non-convex coercive Hamiltonians and study its well-posedness and stability properties. Let us quote that in the latest work, the main results do not require any convexity conditions on the Hamiltonians, contrary to all the previous literature that is based on the control (deterministic) theoretical interpretation of the problem. Among the long list of references on this topic with convex Hamiltonians, we can cite for instance: \cite{control 1}, \cite{control 2}, \cite{control 3}, \cite{control 4}, \cite{control 5}, \cite{control 6}. For recent works on systems of conservative laws posed on junctions, we refer also to \cite{Carda junction 1} and \cite{Carda junction 2} with the references therein. 
In \cite{Lions Souganidis 2}, the authors have studied multi-dimensional junction problems for first and second-order PDE with Kirchhoff-type
Neumann boundary conditions, showing that their generalized viscosity solutions are unique, but still with a vanishing viscosity at the vertex for the second order terms. Finally, let us cite the interesting approach studied in \cite{Lions Souganidis 3}, where it is considered star-shaped tubular domains consisting of a number of non-intersecting, semi-infinite
strips of small thickness that are connected by a central region.  It is shown, that classical regular solutions of uniformly elliptic partial differential equations converge in the thin-domain limit, to the unique solution of a second-order partial differential equation on the network satisfying an effective Kirchhoff-type transmission condition at the junction.
Recall that the key fact in most all of them, is to consider a vanishing viscosity at the vertex $\bf 0$, which is not the case here in the HJB formulation of our stochastic scattering problem of this work.

\subsection{Organization of the paper}
The paper is organized as follows: we introduce in Section \ref{Formulation of the problem and main results} the stochastic optimal scattering control problem and we state the three main results of this contribution, that are proved in Sections \ref{sec: Compactness of the admissible rules and Dynamic Programming Principle}, \ref{sec : theo compa viscosity} and \ref{sec: caracterisation value function}. More precisely, in Section \ref{sec: Compactness of the admissible rules and Dynamic Programming Principle}, we prove the compactness of the admissible rules, and we show the dynamic programming principle in the spirit of \cite{ElKaroui}, adapted to our context. Thereafter, adapting the results obtained in \cite{Ohavi Walsh PDE} to our purpose, we prove a comparison theorem for parabolic backward HJB continuous viscosity solutions posed on star-shaped networks, having the {\it non linear local-time Kirchhoff's boundary condition} at $\bf 0$. In Section \ref{sec: caracterisation value function}, we show that the last system studied in the previous Section \ref{sec : theo compa viscosity} characterizes in a unique way the value function of our problem of stochastic scattering control, which is governed by an optimal diffraction probability measure selected by the own local-time of the spider. In the last Section, \ref{sec control without l} we discuss stochastic control problems for spider diffusion without dependency with respect to local-time variable, for which their well-posedness can be deduced as a consequence of the main results of this contribution.

\section{Formulation of the problem and main results}\label{Formulation of the problem and main results}

In this section, we introduce the main notations for our purpose and we state the main results of this work.

\subsection{Notations and definitions}\label{notations}
Fix $I\geq 2$ an integer. Denote $[I]$ in all this work the following set: 
$$[I]:=\{1   \ldots    I\}.$$
We define the star-shaped network $\mathcal{N}$ with $I$ edges by:
\begin{align*}
\mathcal{N}~~:=~~\{{\bf 0}\}\cup \Big((0,\infty)\times [I]\Big).
\end{align*}
All the points of $\mathcal{N}$ are described by couples $(x,i)\in [0,\infty)\times [I]$, whereas the junction point ${\bf 0}$ is identified with the equivalent class $\{(0,i)~:~i\in [I]\}$. We will also often identify the space $\mathcal{N}$ with a union of $I$ edges
$\mathcal{R}_i=[0,+\infty)$ satisfying $\mathcal{R}_i\cap \mathcal{R}_j=\bf 0$, whenever $(i,j)\in [I]^2$ with $i\neq j$. With these notations, $(x,i)\in {\mathcal{N}}$ is equivalent to asserting that $x\in \mathcal{R}_i$.

We endow naturally ${\mathcal{N}}$ with the following geodesic distance $d^{\mathcal{N}}$ defined by
\begin{equation*}
\forall \Big((x,i),(y,j)\Big)\in \mathcal{N}^2,~~d^\mathcal{N}\Big((x,i),(y,j) \Big)  := \left\{
\begin{array}{ccc}
 |x-y| & \mbox{if }  & i=j\;,\\ 
x+y & \mbox{if }  & {i\neq j},\;
\end{array}\right.
\end{equation*}
so that $\big(\mathcal{N}, d^{\mathcal{N}}\big)$ is a Polish space. We denote by $\mathcal{C}^{\mathcal{N}}[0,T]$ the Polish space of maps defined from $[0,T]$ onto the junction space $\mathcal{N}$ that are continuous w.r.t. the metric $d^{\mathcal{N}}$. The space $\mathcal{C}^{\mathcal{N}}[0,T]$ is naturally endowed with the uniform metric $d^\mathcal{N}_{[0,T]}$ defined by: 
\begin{equation*}
\forall \Big((x,i),(y,j)\Big)\in \Big(\mathcal{C}^{\mathcal{N}}[0,T]\Big)^2,~~d^\mathcal{N}_{[0,T]}:=\sup_{t \in [0,T]}d^\mathcal{N}\Big((x(t),i(t)),(y(t),j(t)) \Big).
\end{equation*}
Together with $\mathcal{C}^{\mathcal{N}}[0,T]$, we introduce
\begin{eqnarray*}
\mathcal{L}[0,T]~~:=~~\Big\{l:[0,T]\to [0,+\infty),\text{ continuous and non decreasing}\Big\}
\end{eqnarray*} 
endowed with the usual uniform distance $|\,.\,|_{(0,T)}$. The modulus of continuity on $\mathcal{C}^{\mathcal{N}}[0,T]$ and $\mathcal{L}[0,T]$ are naturally defined for any $\theta\in (0, T]$ as
\begin{align*}
&\forall X=\big(x,i\big)\in \mathcal{C}^{\mathcal{N}}[0,T],\\
&\hspace{1,3 cm}\omega\big(X,\theta\big)= \sup\Big\{d^{\mathcal{N}}\big((x(s),i(s)),(x(u),i(u))\big) ~\big\vert~(u,s)\in [0,T]^2,~~|u-s|\leq \theta\Big\};\\ 
&\forall f\in {\mathcal{L}}[0,T],\\
&\hspace{1,3 cm}\omega\big(f,\theta\big)= \sup \Big\{|f(u)-f(s)|~\big\vert~(u,s)\in [0,T]^2,~~|u-s|\leq \theta\Big\}.
\end{align*}
We then form the product space 
\begin{align*}
\Phi~~=~~\mathcal{C}^{\mathcal{N}}[0,T]\times \mathcal{L}[0,T]
\end{align*}
 considered as a measurable Polish space equipped with its Borel $\sigma$-algebra $\mathbb{B}(\Phi)$ generated by the open sets relative to the metric $d^{\Phi}:=d^{\mathcal{N}}_{[0,T]} + |\,.\,|_{(0,T)}$. Note to simplify the notations, we will often denote $(x,i,l)$ an element of $\mathcal J\times [0,+\infty)$, rather than $((x,i),l)$, and in the same way $\big(x(s),i(s),l(s)\big)_{s\in[0,T]}\in \Phi$.
 
Given a compact $\mathcal{K}$ of $\R^n$ ($n \in \mathbb{N}^*$), we define:
\begin{eqnarray*}
L^{\infty}_{mc}([0,T]\times \mathcal{K}):=\Big\{f\in L^\infty([0,T]\times \mathcal{K}),~~k\mapsto f(t,k)\in \mathcal{C}(\mathcal{K})~~,\forall t\in [0,T]\Big\}.
\end{eqnarray*}
We denote by $\mathcal{M}_{mc}([0,T]\times \mathcal{K})$ the set consisting of non negative finite measures on $\Big([0,T]\times \mathcal{K},\mathbb{B}([0,T])\otimes \mathbb{B}(\mathcal{K})\Big)$, endowed with the finest topology making continuous the following family of linear forms 
$(\theta_f)_{f\in L^{\infty}_{mc}([0,T]\times \mathcal{K})}$:
\begin{eqnarray*}
\theta_f:\begin{cases}
\mathcal{M}_{mc}([0,T]\times \mathcal{K}))\to \R \\
\nu\mapsto \nu(f)=\displaystyle \int_{[0,T]\times \mathcal{K}} f d\nu
\end{cases}.
\end{eqnarray*}

Let us define now the set of controls (generalized actions) for our problem of stochastic scattering control. \\
Let $(\mathcal{B}_i)_{i\in [I]}$ be $I$ compact sets of $\R$. The set of the generalized actions $\mathcal{U}([0,T]\times \mathcal{B}_i)$ on each edge $\mathcal{R}_i$ is defined by: 
\begin{eqnarray*}
\mathcal{U}([0,T]\times \mathcal{B}_i):=\Big\{\nu \in \mathcal{M}_{mc}([0,T]\times \mathcal{B}_i),~~ \nu^{[0,T]}(dt)=\int_{\mathcal{B}_i}\nu(dt,d\beta_i)=dt\Big\}.
\end{eqnarray*}
Remark that for all $i\in[I]$, $\mathcal{U}([0,T]\times \mathcal{B}_i)$ are compact for the weak topology\\
$*\sigma\Big(L^\infty_{mc}([0,T]\times \mathcal{B}_i)^{'},L^\infty_{mc}([0,T]\times \mathcal{B}_i)\Big)$.\\
Given $\mathcal{O}$ a compact set of $\R^I$ , we define $\mathcal{V}([0,T]\times \mathcal{O})$ the set of the generalized actions at the junction point $\bf 0$ by:
\begin{eqnarray*}
&\mathcal{V}([0,T]\times \mathcal{O}):=\Big\{\nu \in \mathcal{M}_{mc}([0,T]\times \mathcal{O}),~~\exists \ell_\nu\in \mathcal{L}([0,T]),~\text{s.t}~~\nu^{[0,T]}(dt)=d\ell_\nu(t)\Big\},\\
&\text{where }\nu^{[0,T]}(dt)=\displaystyle\int_{\mathcal{O}}\nu(dt,d\vartheta).
\end{eqnarray*}
Remark that as a consequence of the disintegration Theorem of a measure, (see for instance \cite{mesure inte}), we have for all $\nu \in \mathcal{V}([0,T]\times \mathcal{O})$: $$\nu(dt,d\vartheta)=d\ell_\nu(t)\nu_{t}(d\vartheta),$$ where $\nu_.$ is a measurable kernel of mass $1$ on $(\mathcal{O},\mathbb{B}(\mathcal{O}))$.

\textbf{Fix now in all of this work, $I$ compact sets $(\mathcal{B}_i)_{i\in[I]}$ of $\R$ and $\mathcal{O}$ a compact set of $\R^I$.}

The canonical space involved in the martingale formulation for our purpose is defined by:
$$
\Phi~~=~~\mathcal{C}^{\mathcal{N}}([0,T])\times\mathcal{L}([0,T])\times \Big(\prod_{i=1}^I \mathcal{U}([0,T]\times \mathcal{B}_i)\Big)\times \mathcal{V}([0,T]\times \mathcal{O}),
$$ 
and we will naturally denote by $\mathbb{B}(\Phi)$ its $\sigma$ Borel algebra.

The canonical \textbf{continuous} process defined on $\big(\Phi,\mathbb{B}(\Phi)\big)$ is given by:
\begin{eqnarray} \label{eq : def processus cano}
& X:\begin{cases}
[0,T]\times\Phi \to \mathcal{N}\times[0,+\infty)\times\Big(\displaystyle\prod_{i=1}^I \mathcal{M}_{mc}([0,T]\times \mathcal{B}_i)\Big)\times \mathcal{M}_{mc}([0,T]\times \mathcal{O})\\
\big(s,\omega\big)\mapsto \Big(x(s),i(s),l(s),\nu_1(s) \ldots   \nu_I(s),\nu_0(s)\Big),\\
\end{cases}
\\
&\nonumber \text{where:}~~\omega =\Big(\big(x(s),i(s),l(s)\big)_{s\in[0,T]},\nu_1   \ldots   \nu_I,\nu_0\Big),\\&\nonumber \nu_i(s)(dt,d\beta_i)=\mathbf{1}_{[0,s]}(t)\nu_i(dt,d\beta_i),~~\forall s\in [0,T],~~\forall i \in [I],\\
&\nonumber \nu_0(s)(dt,d\theta_1, \ldots ,d\theta_I)=\mathbf{1}_{[0,s]}(t)\nu_0(dt,d\theta_1,\ldots,d\theta_I),~~\forall s\in [0,T].
\end{eqnarray}
We denote by $(\Psi_t=\sigma(X(s),~0\leq s \leq t))_{0\leq t\leq T}$ the right continuous filtration generated by the process $X(\cdot)$ on $\big(\Phi, \mathbb{B}(\Phi)\big)$. 
The set of the probability measures on the measurable space $(\Phi,\mathbb{B}(\Phi))$ is denoted by $\mathcal{P}(\Phi,\mathbb{B}(\Phi))$. We endow this space with the topology of weak convergence of probability measures. Since $\Phi$ is a Polish space, this topology is induced by the Prokhorov metric $\pi$ on $\mathcal{P}(\Phi,\mathbb{B}(\Phi))$. Recall that the Prokhorov distance $\pi(P,Q)$ between two probability measures on $(\Phi,\mathbb{B}(\Phi))$ is defined as
\begin{align*}
&\forall P,Q \in \mathcal{P}(\Phi,\mathbb{B}(\Phi)),\\
&\hspace{0.5 cm}\pi(P,Q)=\inf\left \{\varepsilon>0~\vert~\forall A\in \mathbb{B}(\Phi),\;\; P(A)\leq Q(A^{\varepsilon}) + \varepsilon \text{ and }Q(A)\leq P(A^{\varepsilon}) + \varepsilon\right \}
\end{align*}
where $A^{\varepsilon}$ denotes the $\varepsilon$-dilation of the set $A\in \Phi$ ($A^{\varepsilon}=\{\tilde{X}\in \Phi,~d^{\Phi}(\tilde{X}, A)\leq \varepsilon\}$).
We introduce the set $\mathcal{P}([I])\subset [0,1]^I$   \begin{align*}
    \mathcal{P}([I]):=\Big\{(\mathcal{S}_{i})\in[0,1]^{I}~\big\vert~ \sum_{i=1}^{I}\mathcal{S}_i=1\Big\}
\end{align*}
the simplex set giving all probability measures on $[I]$.

\subsection{Martingale problem and viscosity solutions:}
In this subsection, we define our stochastic scattering control problem in the weak sense, using a martingale formulation. We also introduce the value function of the  stochastic control problem. Finally, we give the necessary tools in the viscosity theory adapted to our purpose, and we introduce the non linear HJB system; that will drive the optimal scattering process at $\bf 0$.\\
Recall that $(\mathcal{B}_i)_{i\in[I]}$ is a collection of $I$ compact sets of $\R$, whereas  $\mathcal{O}$ is a compact set of $\R^I$. 

In the sequel, we define the following domain:
$$\mathcal{D}_T=[0,T]\times \mathcal{N}\times [0,+\infty).$$

Let us introduce $\mathcal{C}^{1,2,1}_b(\mathcal{D}_T)$ the class of continuous function defined on $\mathcal{D}_T$ with regularity $\mathcal{C}^{1,2,1}_b([0,T]\times [0,\infty)^2)$ on each edge, namely
\begin{align*}
&\mathcal{C}^{1,2,1}_b(\mathcal{D}_T):=\Big\{f:\mathcal{J}_T\to\R,~~(t,(x,i),l)\mapsto f_i(t,x,l)~\Big\vert~\forall i\in [I],\\
&f_i:[0,T]\times {J}_i\times[0+\infty)\to\R,\,(t,x,l)\mapsto f_i(t,x,l)\in \mathcal{C}^{1,2,1}_b([0,T]\times J_i\times[0+\infty)),\\
&\hspace{3 cm}~ \forall (t,(i,j),l)\in [0,T]\times [I]^{2}\times[0+\infty), \,f_i(t,0,l)=f_j(t,0,l)\Big\}.
\end{align*}

Now, we introduce the following data:
$$\textbf{Data:~} (\mathcal{D})~~\begin{cases}
\Big(\sigma_i \in \mathcal{C}\big([0,T]\times[0,+\infty)^2\times \mathcal{B}_i,\R\big)\Big)_{i\in[I]}\\
\Big(b_i \in \mathcal{C}\big([0,T]\times[0,+\infty)^2\times\mathcal{B}_i,\R\big)\Big)_{i\in[I]}\\
\Big(h_i \in \mathcal{C}\big([0,T]\times[0,+\infty)^2\times \mathcal{B}_i,\R\big)\Big)_{i\in[I]}\\
\mathcal{S}=\big(\mathcal{S}_i\big)_{i\in[I]}\in \mathcal{C}\big([0,T]\times[0+\infty)\times\mathcal{O},\mathcal{P}[I]\big)\\
h_0\in \mathcal{C}\big([0,T]\times[0+\infty)\times\mathcal{O},\R\big),\\
\Big(g_i \in \mathcal{C}\big([0,+\infty)^2,\R\big)\Big)_{i\in[I]}
\end{cases}.
$$
We assume that the data $(\mathcal{D})$ satisfy the following assumptions: (where $({\bf A})$ stands for the spider spinning measure $\mathcal{S}$, $({\bf E})$ for ellipticity, and $({\bf R})$ for Lipschitz regularity uniformly with respect to the control variables):
\newpage
$$\textbf{Assumption } (\mathcal{H})$$
\begin{align*}
&({\bf A})~~\exists\;\ds \underline{\zeta}>0,~~\forall i\in[I], ~~\forall (t,l,\vartheta)\in [0,T]\times[0,+\infty)\times\mathcal{O},~~\mathcal{S}_i(t,l,\vartheta)~~\ge~~\underline{\zeta}.\\
&({\bf E})~~\exists\,\,\underline{\sigma} >0,~~\forall i\in[I], ~~\forall (t,x,l,\beta_i)\in [0,T]\times[0,+\infty)^2\times \mathcal{B}_i,~~\sigma_i(t,x,l,\beta_i)~~\ge~~\underline{\sigma}.\\
&({\bf R})~~\exists\,\,(|b|,|h|,\overline{\zeta},\overline{\sigma})\in (0,+\infty)^4,~~\forall i\in[I],\\
&({\bf R} - i)\displaystyle\sup_{t,x,l,\beta_i}|b_i(t,x,l,\beta_i)|+\sup_{t,l,\beta_i}\sup_{(x,y),\;x\neq y} \frac{|b_i(t,x,l,\beta_i)-b_i(t,y,l,\beta_i)|}{|x-y|}\\
&\displaystyle\hspace{3,5cm}+\sup_{x,l,\beta_i}\sup_{(t,s),\;t\neq s} \frac{|b_i(t,x,l,\beta_i)-b_i(s,x,l,\beta_i)|}{|t-s|}\\
&\displaystyle\hspace{3,5cm}+\sup_{t,x,\beta_i}\sup_{(l,l'),\;l\neq l'} \frac{|b_i(t,x,l,\beta_i)-b_i(t,x,l',\beta_i)|}{|l-l'|}~~\leq~~|b|,\\
&({\bf R} - ii)\displaystyle\sup_{t,x,l,\beta_i}|\sigma_i(t,x,l,\beta_i)|+\sup_{t,l,\beta_i}\sup_{(x,y),\;x\neq y} \frac{|\sigma_i(t,x,l,\beta_i)-\sigma_i(t,y,l,\beta_i)|}{|x-y|}\\
&\displaystyle\hspace{3,5cm}+\sup_{x,l,\beta_i}\sup_{(t,s),\;t\neq s} \frac{|\sigma_i(t,x,l,\beta_i)-\sigma_i(s,x,l,\beta_i)|}{|t-s|}\\
&\displaystyle\hspace{3,5cm}+\sup_{t,x,\beta_i}\sup_{(l,l'),\;l\neq l'} \frac{|\sigma_i(t,x,l,\beta_i)-\sigma_i(t,x,l',\beta_i)|}{|l-l'|}~~\leq~~\overline{\sigma},\\
&({\bf R} - iii) \displaystyle\sup_{t,x,l,\beta_i}|h_i(t,x,l,\beta_i)|+\sup_{t,l,\beta_i}\sup_{(x,y),\;x\neq y} \frac{|h_i(t,x,l,\beta_i)-h_i(t,y,l,\beta_i)|}{|x-y|}\\
&\displaystyle\hspace{3,5cm}+\sup_{x,l,\beta_i}\sup_{(t,s),\;t\neq s} \frac{|h_i(t,x,l,\beta_i)-h_i(s,x,l,\beta_i)|}{|t-s|}\\
&\displaystyle\hspace{3,5cm}+\sup_{t,x,\beta_i}\sup_{(l,l'),\;l\neq l'} \frac{|b_i(t,x,l,\beta_i)-b_i(t,x,l',\beta_i)|}{|l-l'|}~~\leq~~|h|,\\
&({\bf R}- iv)~\displaystyle\sup_{t,l,\vartheta}|\mathcal{S}_i(t,l,\vartheta)|+\sup_{l,\vartheta}\sup_{(t,s),\;t\neq s} \frac{|\mathcal{S}_i(t,l,\vartheta)-\mathcal{S}_i(s,l,\vartheta)|}{|t-s|}+\\
&\hspace{3,5cm}+\sup_{t,\vartheta}\sup_{(l,l'),\;l\neq l'} \frac{|\mathcal{S}_i(t,l,\vartheta)-\mathcal{S}_i(t,l',\vartheta)|}{|l-l'|}~~\leq ~~\overline{\zeta},\\
&({\bf R}- v)~\displaystyle\sup_{t,l,\vartheta}|h_0(t,l,\vartheta)|+\sup_{l,\vartheta}\sup_{(t,s),\;t\neq s} \frac{|h_0(t,l,\vartheta)-h_0(s,l,\vartheta)|}{|t-s|}+\\
&\hspace{3,5cm}+\sup_{t,\vartheta}\sup_{(l,l'),\;l\neq l'} \frac{|h_0(t,l,\vartheta)-h_0(t,l',\vartheta)|}{|l-l'|}~~\leq ~~|h|.
\end{align*}

Fix now $(t,x,i,l)\in [0,T]\times \mathcal{N}\times [0,+\infty)$. 
\textbf{We define the set of admissible rules $\mathcal{A}\big(t,x,i,l\big)$}, as the set of all the 'possible controlled' - probability measures $\P^{t,x,i,l}_{\beta,\vartheta}\in \mathcal{P}\Big(\Phi,\mathbb{B}(\Phi)\Big)$, such that the following conditions $(\mathcal{S})$ are satisfied:
$$ \textbf{ Conditions } (\mathcal{S})$$
(i) For each $u\leq t$: $$\Big(x(u),i(u),l(u),(\nu_i(u))_{i\in [I]},\nu_0(u)\Big)=\Big(x,i,l,(\nu_i(t))_{i\in [I]},\nu_0(t)\Big),~~\P^{t,x,i,l}_{\beta,\vartheta}~~ \text{a.s.}$$ 
(ii) For each $s\ge t$, 
\begin{eqnarray*}
&\nu^{[0,T]}_0(dt)=\displaystyle\int_{\mathcal{O}}\nu_0(dt,d\vartheta)=dl(t)=d\ell_{\nu_0}(t),\\
&\displaystyle\int_{t}^{s} \int_{\mathcal{O}}\mathbf{1}_{\{x(u)>0\}}\nu_0(s)(du,d\vartheta)=\displaystyle\int_{t}^{s}\mathbf{1}_{\{x(u)>0\}}dl(s)=0,~~\P^{t,x,i,l}_{\beta,\vartheta} \text{ a.s.}
\end{eqnarray*}
(iii) For any $f\in\mathcal{C}^{1,2,1}_{b} \big(\mathcal{D}_T\big)$, the following process:
\begin{eqnarray*}
&\Big(~~f_{i(s)}(s,x(s),\ell(s))- f_{i}(t,x,l)
-\\
&\displaystyle\int_{t}^{s}\int_{\mathcal{B}_{i(u)}}\Big[\partial_tf_{i(u)}(u,x(u),l(u))+
b_{i(u)}(u,x(u),l(u),\beta_{i(u)})\partial_xf_{i(u)}(u,x(u),l(u))+\\
&\displaystyle \frac{\sigma_{i(u)}^2(u,x(u),l(u),\beta_{i(u)})}{2}\partial_x^2f_{i(u)}(u,x(u),l(u))\Big]\nu_i(s)(du,d\beta_{i(u)})-\ds\int_{t}^{s}\Big[\partial_lf(u,0,l(u))+\\
&\displaystyle\sum_{i=1}^{I}\int_{\mathcal{O}}\mathcal{S}_i(u,l(u),\vartheta)\partial_{x}f_{i}(u,0,l(u))\Big]\nu_0(s)(du,d\vartheta)~~\Big)_{t\leq s\leq T},
\end{eqnarray*} 
is a $(\Psi)_{t\leq s\leq T}$ continuous martingale under the probability measure $\P^{t,x,i,l}_{\beta,\vartheta}$, after time $t$.
\begin{Remark}
As a straight consequence of the main result obtained in \cite{Spider} together with assumption $(\mathcal{H})$, note that $\mathcal{A}\big(t,x,i,l\big)$ is non empty.    
\end{Remark}

Next, we introduce the value function of the stochastic scattering control problem.
\begin{Definition}\label{def fonction valeur}
Assume assumption $(\mathcal{H})$. The value function $u$ of the stochastic scattering control problem is defined by:
\begin{eqnarray}\label{eq : value function}
&\nonumber u:=\\
&\begin{cases}
[0,T]\times \mathcal{N}\times[0,+\infty)\to \R \\
\Big(t,x,i,l\Big)\mapsto \sup \Big\{~~\mathbb{E}^{\P^{t,x,i,l}_{\beta,\vartheta}}\Big[~~\displaystyle\displaystyle\int_{t}^{T}\int_{\mathcal{B}_{i(u)}}h_{i(u)}(u,x(u),l(u),\beta_{i(u)})\nu_{i(u)}(T)(du,d\beta_{i(u)})\\
\hspace{3.2 cm}+\displaystyle \int_{t}^{T}\int_{\mathcal{O}}h_0(u,l(u),\vartheta)\nu_0(T)(du,d\vartheta)~+~g_{i(T)}\big(x(T),\ell(T)\big)~~\Big],\\
\hspace{9.5 cm}\P^{t,x,i,l}_{\beta,\vartheta}\in\mathcal{A}\big(t,x,i,l\big)~~\Big\} 
\end{cases}
\end{eqnarray}
where $h_0$ is the cost of the agent at the junction point $\bf 0$, the $h_i$ correspond to the costs on each ray $\mathcal{R}_i$ $(i\in [I])$, and $g$ is the terminal condition. 
\end{Definition}
We introduce now the HJB system that will characterize the problem of control; called {\it Walsh's spider HJB system with non linear local-time Kirchhoff's boundary condition} in \cite{Ohavi Walsh PDE}:
\begin{eqnarray}\label{eq PDE Walsh}
&\begin{cases}
\textbf{HJB equation parameterized by the local-time on each ray}:\\
\partial_t u_i(t,x,l)+\underset{\beta_i\in \mathcal{B}_i}{\sup}\Big\{\ds \frac{1}{2}\sigma_i^2(t,x,l,\beta_i)\partial^2_xu_i(t,x,l)+\\
\hspace{1 cm}b_i(t,x,l,\beta_i)\partial_xu_i(t,x,l)+h_i(t,x,l,\beta_i)\Big\}=0,~~(t,x,l)\in(0,T)\times \R_+^2,\\
\textbf{Non linear local-time Kirchhoff's boundary transmission at } \bf 0: \\
\partial_lu(t,0,l)+\underset{ \vartheta \in \mathcal{O}}{\sup} \Big\{\ds\sum_{i=1}^I\mathcal{S}_i(t,l,\vartheta)\partial_xu_i(t,0,l)\\
\hspace{3 cm}+h_0(t,l,\vartheta)\Big\}=0,~~(t,l)\in(0,T)\times\R_+,~~\\
\textbf{Terminal condition}:\\
u_i(T,x,l)=g_i(x,l),~~\forall (x,l)\in\R_+^2.~~\forall i\in[I],\\
\textbf{Continuity condition at } \bf 0:\\
\forall (i,j)\in[I]^2,~~\forall (t,l)\in \R_+,~~u_i(t,0,l)=u_j(t,0,l).
\end{cases}
\end{eqnarray}
We end this section by giving the definition of continuous viscosity super and sub solutions that belong to $\mathcal{C}\big(\mathcal{D}_T\big)$ (recall that $\mathcal{D}_T=[0,T]\times \mathcal{N}\times [0,+\infty)$), defined by:
\begin{align*}
&\mathcal{C}\big(\mathcal{D}_T\big):=\Big\{f:[0,T]\times \mathcal{N}\times [0,+\infty),~~(t,x,i,l)\mapsto f_i(t,x,l)~\Big|\\
&~\forall i\in [I],~f_i:[0,T]\times \mathcal{N}\times [0,+\infty)\to\R,\,(t,x,l)\mapsto f_i(t,x,l)\in \mathcal{C}^{0}_b([0,T]\times \mathcal{N}\times [0,+\infty)),\\
&~\forall (i,j,t,l)\in [I]^{2}\times [0,T]\times [0,+\infty), ~f_i(t,0,l)=f_j(t,0,l)=f(t,0,l) \Big\}.
\end{align*} 
for the {\it Walsh's spider backward HJB system}, given in \eqref{eq PDE Walsh}.
In order to remain consistent with the results obtained in \cite{Martinez-Ohavi EDP}, more precisely with the class of regularity of the solutions in the linear framework (see Introduction in \cite{Martinez-Ohavi EDP}), we introduce the following space of test functions for continuous viscosity solutions of system \eqref{eq PDE Walsh}:
\begin{align*}
&\mathcal{C}^{1,2,0}_{{\bf 0},1}\big(\mathcal{D}_T\big):=\Big\{f:\mathcal{D}_T,~~(t,x,i,l)\mapsto f_i(t,x,l)~\Big\vert\\
&\hspace{2 cm}~\forall i\in [I], \;f_i:[0,T]\times [0,+\infty)^2\to\R,\\&\hspace{2 cm}~(t,x,l)\mapsto f_i(t,x,l)\in \mathcal{C}^{0,1,0}_b([0,T]\times (0,+\infty)^2)\cap \mathcal{C}^{1,2,0}_b((0,T)\times (0,+\infty)^2),\\
&\hspace{2 cm}~ \forall (i,j,t,l)\in [I]^{2}\times [0,T]\times[0,+\infty), \,f_i(0,l)=f_j(0,l)=f(0,l),
\\&\hspace{2 cm}~f(\cdot,0,\cdot)\in \mathcal{C}^{0,1}_b\big([0,T]\times [0,+\infty)\big)\Big\}.
\end{align*}
\begin{Definition}\label{def viscosity}
Let $u\in \mathcal{C}\big(\mathcal{D}_T\big)$.\\
a) We say that $u$ is a continuous viscosity super solution of the HJB system \eqref{eq PDE Walsh}, if for all test function $\phi \in \mathcal{C}^{1,2,0}_{{\bf 0},1}\big(\mathcal{D}_T\big)$ and for all local minimum point\\
$(t_\star,x_\star,i_\star,l_\star)\in [0,T]\times [I] \times [0,+\infty)^2$ of $u-\phi$,  with $(u-\phi)_{i_\star}(t_\star, x_\star,l_\star)=0$, we have:
\begin{align*}
\begin{cases}
\partial_t\phi_{i_\star}(t_\star,x_\star,l_\star)+\underset{\beta_{i_\star}\in \mathcal{B}_{i_\star}}{\sup}\Big\{\ds \frac{1}{2}\sigma_{i_\star}(x_\star,l_\star,\beta_i)^2\partial^2_x\phi_{i_\star}(t_\star,x_\star,l_\star)+\\
b_{i_\star}(x_\star,l_\star,\beta_i)\partial_x\phi_{i_\star}(t_\star,x_\star,l_\star)+h_{i_\star}(x_\star,l_\star,\beta_i)\Big\}\leq 0,~\text{if}~(t_\star,x_\star,l_\star)\in(0,T)\times (0,+\infty)^2,\\
\partial_l\phi(t_\star,0,l_\star)+
\underset{\vartheta \in \mathcal{O}}{\sup} \Big\{\ds \sum_{i=1}^I\mathcal{S}_i(t_\star,l_\star,\vartheta)\partial_x\phi_i(t_\star,0,l_\star)\\
+h_0(t_\star,l_\star,\vartheta)\Big\} \leq 0,~\text{if}~x_\star=0.~~(t_\star,l_\star)\in(0,T)\times (0,+\infty),
\end{cases}.
\end{align*}
b) We say that $u$ is a continuous viscosity sub solution of the HJB system \eqref{eq PDE Walsh}, if for all test function $\phi \in \mathcal{C}^{1,2,0}_{{\bf 0},1}\big(\mathcal{D}_T\big)$ and for all local maximum point\\
$(t_\star,x_\star,i_\star,l_\star)\in [0,T]\times [I] \times [0,+\infty)^2$ of $u-\phi$, with $(u-\phi)_{i_\star}(t_\star, x_\star,l_\star)=0$, we have:
\begin{align*}
\begin{cases}
\partial_t\phi_{i_\star}(t_\star,x_\star,l_\star)+\underset{\beta_{i_\star}\in \mathcal{B}_{i_\star}}{\sup}\Big\{\ds \frac{1}{2}\sigma_{i_\star}(x_\star,l_\star,\beta_i)^2\partial^2_x\phi_{i_\star}(t_\star,x_\star,l_\star)+\\
b_{i_\star}(x_\star,l_\star,\beta_i)\partial_x\phi_{i_\star}(t_\star,x_\star,l_\star)+h_{i_\star}(x_\star,l_\star,\beta_i)\Big\}\ge 0,~\text{if}~(t_\star,x_\star,l_\star)\in(0,T)\times (0,+\infty)^2,\\
\partial_l\phi(t_\star,0,l_\star)+
\underset{\vartheta \in \mathcal{O}}{\sup} \Big\{\ds \sum_{i=1}^I\mathcal{S}_i(t_\star,l_\star,\vartheta)\partial_x\phi_i(t_\star,0,l_\star)\\
+h_0(t_\star,l_\star,\vartheta)\Big\} \ge 0,~\text{if}~x_\star=0.~~(t_\star,l_\star)\in(0,T)\times (0,+\infty),
\end{cases}.
\end{align*}
c) We say that $u$ is a continuous viscosity solution of
the HJB system \eqref{eq PDE Walsh}, if it is both a continuous viscosity super and sub solution of the HJB system \eqref{eq PDE Walsh}.
\end{Definition}

\subsection{Main results}\label{subs Main results}
We conclude this section by stating the main results of this work, that are:\\
-the dynamic programming principle satisfied by the value function $u$, Theorem \ref{th: PPD},\\
-the comparison theorem, Theorem \ref{th compa visco} for continuous viscosity solution of our Walsh's spider HJB system \eqref{eq PDE Walsh}, having {\it non linear local-time Kirchhoff's boundary condition at $\bf 0$}.\\
-the unique characterization of the value function $u$ given \eqref{eq : value function} in terms of the HJB system \eqref{eq PDE Walsh}.
\begin{Theorem}\label{th: PPD}
\textbf{Dynamic programming principle equation:} Assume assumption $(\mathcal{H})$.  For any $(\Psi_s)_{t\leq s \leq T}$ stopping time $\tau$, the value function $u$ defined in Definition \ref{def fonction valeur}, satisfies for all $(i,t,x,l) \in[I]\times[0,T]\times[0,+\infty)^2$:
\begin{align}\label{eq: ppd}
&\nonumber u_i(t,x,l)=\sup \Big\{~~\mathbb{E}^{\P^{t,x,i,l}_{\beta,\vartheta}}\Big[~~\displaystyle\int_{t}^{\tau}\int_{\mathcal{B}_{i(u)}}h_{i(u)}(u,x(u),l(u),\beta_{i(u)})\nu_{i(u)}(T)(du,d\beta_{i(u)})+\\
&\nonumber \int_{t}^{\tau}\int_{\mathcal{O}}h_0(u,l(u),\vartheta)\nu_0(T)(du,d\vartheta)+u_{i(\tau)}\big(\tau,x(\tau),l(\tau)\big)~~\Big],~~ \P^{t,x,i,l}_{\beta,\vartheta}\in\mathcal{A}\big(t,x,i,l\big)~~\Big\}.
\end{align}
\end{Theorem}
As a consequence of the main result obtained in the pioneered work \cite{Ohavi Walsh PDE} (see Theorem 2.2), adapted to our framework for the stochastic scattering control problem, we have the following comparison theorem:

\begin{Theorem}\label{th compa visco}\textbf{Comparison Theorem:} Assume assumption $(\mathcal{H})$. Let $v\in \mathcal{C}\big(\mathcal{D}_T\big)$ a continuous viscosity sub solution and $u\in \mathcal{C}\big(\mathcal{D}_T\big)$ a continuous viscosity super solution of the Walsh's spider backward HJB system - given in \eqref{eq PDE Walsh}, satisfying the following backward boundary conditions:
\begin{align*}
&\forall i\in [I],~~\forall (x,l)\in[0,+\infty)^2,~~u_i(T,x,l)\ge v_i(T,x,l),
\end{align*}
Assume that $u$ and $v$ are uniformly bounded in the whole domain $\mathcal{D}_T$. Then we have:
$$\forall (t,i,x,l)\in [0,T]\times i\in [I]\times [0,+\infty)^2,~~u_i(t,x,l)\ge v_i(t,x,l). $$
\end{Theorem}
\begin{Remark}\label{rm: non vide}
As it is classical in the proofs of comparison theorems in unbounded domains, one can also impose that $u$ and $v$ satisfy a linear growth, rather than be uniformly bounded. To lighten the article, we have preferred the last assumption. 
\end{Remark}
The third and last main result of this work is the following:
\begin{Theorem}\label{th unique caracterization for u}\textbf{Unique characterization of the value function of the stochastic scattering problem, with optimal diffraction probability measure at the vertex selected from the own local time:}\\
Assume assumption $(\mathcal{H})$. The value function $u$ defined in Definition \ref{def fonction valeur} is the unique continuous viscosity solution of the Walsh's spider backward HJB system - given in \eqref{eq PDE Walsh}.
\end{Theorem}
Let us conclude by an important remark.
\begin{Remark}
In the next Section, we will prove that the value function $u$ defined in \eqref{eq : value function}, indeed attains its maximum (Theorem \ref{th: compact}). In other terms, there exists an optimal control $\overline{\P}^{t,x,i,l}_{\overline{\beta},\overline{\vartheta}}$ such that for all $(i,t,x,l)\in[I]\times[0,T]\times[0,+\infty)^2$:
\begin{align}
    &\nonumber u_i(t,x,l)=\mathbb{E}^{\overline{\P}^{t,x,i,l}_{\overline{\beta},\overline{\vartheta}}}\Big[~~\displaystyle\int_{t}^{T}\int_{\mathcal{B}_{i(u)}}h_{i(u)}(u,x(u),l(u),\overline{\beta}_{i(u)})\overline{\nu}_{i(u)}(T)(du,d\overline{\beta}_{i(u)})+\\
&\nonumber \int_{t}^{T}\int_{\mathcal{O}}h_0(u,l(u),\overline{\vartheta})\overline{\nu}_0(T)(du,d\overline{\vartheta})+g_{i(T)}\big(T,x(T),l(T)\big)~~\Big].
\end{align}
Moreover, as a consequence of assumption $(\mathcal{H})$ and the main Theorem contained in \cite{Spider}, it follows that $\overline{\P}^{t,x,i,l}_{\overline{\beta},\overline{\vartheta}}$ is unique (in the sense of distribution). On the other side, the compactness of the set $(\mathcal{B}_i)_{i\in[I]}$ and $\mathcal{O}$, together with the continuity of the data insures existence of optimal feed back controls:
\begin{align*} \beta_i^*(t,x,l,p,S),~~\mathcal{S}_i^*(t,x,l,p_1,\ldots,p_I),\\(i,t,x,l,p.S)\in[I]\times[0,T]\times[0,+\infty)^4,~~(p_1,\ldots,p_I)\in[0,+\infty)^4,
\end{align*}
on each ray $\mathcal{R}_i$ and at the junction point $\bf 0$. As the value function $u$ is not regular enough, it is well know from the theory of stochastic control, that one has to carefully build the optimal feed back control as a function of the upper/sub jets of $u$. Even if we believe that classical arguments can be used on each ray, the case of the junction point $\bf 0$ because its novelty and complexity has to be rigorously investigated. So as not to overload this contribution already quite long, as explained in Introduction, we prefer to postpone this technical point.
\end{Remark}
\section{Compactness of the admissible rules and Dynamic Programming Principle}\label{sec: Compactness of the admissible rules and Dynamic Programming Principle}
In this section, we will prove the compactness of the set of admissible rules $\mathcal{A}\big(t,x,i,l\big)$, and the dynamic programming principle. 
\subsection{A criterion of compactness of the admissible rules at the junction point $\bf 0$.} \label{subs A criterium of compactness of the admissible rules at the junction point}
We start by giving a criterion of compactness for the set of generalized actions $\mathcal{V}([0,T]\times \mathcal{O})$ at the junction point $\bf 0$, that will be used in the proof of the compactness of the admissible rules $\mathcal{A}\big(t,x,i,l\big)$.
\begin{Theorem}\label{th : compactV}
Let $\mathcal{Z}$ be a subset of $\mathcal{V}([0,T]\times \mathcal{O})$. Assume that there exists a constant $C>0$, and a modulus of continuity $m\in\mathcal{C}(\mathbb{R}_{+},\mathbb{R})$, with $m(0)=0$, such that:
\begin{eqnarray*}
~&\forall \nu \in \mathcal{Z},~~ \ell_\nu(T)~~\leq~~C,
\\ 
&\forall \nu \in \mathcal{Z},~~\forall (t,s)\in[0,T],~~ |\ell_\nu(t)-\ell_\nu(s)|~~\leq~~m(|t-s|),
\end{eqnarray*}
then $\mathcal{Z}$ is compact for the weak topology $*\sigma\Big(L^\infty_{mc}([0,T]\times \mathcal{O})^{'},L^\infty_{mc}([0,T]\times \mathcal{O})\Big)$.
\end{Theorem}
\begin{proof}
Recall that the $\sigma$ Borel algebra $\mathbb{B}([0,T])$ of $[0,T]$ is countably generated, hence it follows from Proposition \ref{pr : metri topo faible} that $\mathcal{M}_{mc}([0,T]\times \mathcal{O})$ is metrizable. Therefore $\mathcal{Z}$ is metrizable and the compactness can be proved sequentially.\\
Let $(\nu_n)$ be a sequence of $\mathcal{V}$. By definition, there exists a sequence $(\ell_{\nu_n})$ of $\mathcal{L}([0,T])$, such that: 
\begin{eqnarray*}
\nu_n^{[0,T]}(dt)=\int_{\mathcal{O}}\nu_n(dt,d\vartheta)=\ell_{\nu_n}(dt).
\end{eqnarray*}
Using the assumptions satisfied by the sequence $(\ell_{\nu_n})$ stated in the Proposition; Ascoli's Theorem implies that $(\ell_{\nu_n})$ converges uniformly up to a subsequence to $\ell\in \mathcal{C}[0,T]$, and since $\mathcal{L}([0,T])$ is closed in $\mathcal{C}[0,T]$ for the uniform convergence, we have $\ell\in \mathcal{L}([0,T])$.\\
Let us now show that $\mathcal{Z}$ is relatively compact for $*\sigma\Big(L^\infty_{mc}([0,T]\times \mathcal{O})^{'},L^\infty_{mc}([0,T]\times
\mathcal{O})\Big)$, and for this purpose we will use Theorem \ref{th : compa rela}. It is indeed enough to show that $(\nu_n^{[0,T]})$ and (resp. $(\nu_n^{\mathcal{O}}=\displaystyle\int_{[0,T]}\nu_n(dt,d\vartheta)$)) are relatively compact in $\mathcal{M}_{m}([0,T])$ (resp. $\mathcal{M}_c(\mathcal{O})$), for the weak topologies
$*\sigma\Big(L^{\infty,1}([0,T])^{'},L^{\infty,1}([0,T])\Big)$, \Big(resp. $*\sigma\Big(\mathcal{C}(\mathcal{O})^{'},\mathcal{C}(\mathcal{O})\Big)$\Big). Recall that: 
\begin{eqnarray*}
L^{\infty,1}([0,T])~~:=~~\Big\{~~f\in L^{\infty}([0,T]),~~\exists B\in \mathbb{B}([0,T]),~~f(t)=\mathbf{1}_{B}(t)~~\Big\},
\end{eqnarray*}
and $\mathcal{M}_m([0,T])$, (resp.$\mathcal{M}_c(\mathcal{O})$) are the set of finite positive finite measures on $[0,T]$ (resp. $\mathcal{O}$), endowed with the finest topology making continuous the following family of linear forms $(\theta_f)_{f\in L^{\infty}([0,T])}$: 
\begin{eqnarray*}
\theta_f:=\begin{cases}
\mathcal{M}_m([0,T])\to \R \\
\nu\mapsto \nu(f)=\displaystyle \int_{[0,T]} f d\nu
\end{cases}.
\end{eqnarray*}
\Big(resp. $(\theta_f)_{f\in \mathcal{C}(\mathcal{O})}$
\begin{eqnarray*}
\theta_f:=\begin{cases}
\mathcal{M}_c(\mathcal{O})\to \R \\
\nu\mapsto \nu(f)=\displaystyle \int_{\mathcal{O}} f d\nu
\end{cases}\Big).
\end{eqnarray*}
Since $(\ell_n)$ converges uniformly to $\ell$ up to a subsequence $n_k$, it is easy to get that for any $f\in L^{\infty,1}([0,T])$:
\begin{eqnarray*}
\displaystyle\int_{[0,T]} f(t)\ell_{\nu_{n_k}}(dt)~~\xrightarrow[n_k]{k\to +\infty}~~\displaystyle\int_{[0,T]}  f(t)\ell(dt),
\end{eqnarray*}
in other terms $\nu_{n_k}^{[0,T]}(dt)\overset{*}{\rightharpoonup}\ell(dt)$ for $*\sigma\Big(L^{\infty,1}([0,T])^{'},L^{\infty,1}([0,T])\Big)$.\\
On the other hand, we have: 
\begin{eqnarray*}
\|\nu_{n}^{\mathcal{O}}\|_{\mathcal{C}(\mathcal{O})^{'}}~~=\sup_{f\in \mathcal{C}(\mathcal{O}), \|f\|\leq 1 }~~\Big|\int_{[0,T]\times \mathcal{O}}f(t)\nu_{n}(dt,d\vartheta)\Big|~~\leq~~\ell_n(T)~~\leq~~C,\end{eqnarray*}
and this implies that $(\nu_{n}^{\mathcal{O}})$ is relatively compact for the weak topology $*\sigma(\mathcal{C}(\mathcal{O})^{'},\mathcal{C}(\mathcal{O}))$.
We deduce finally using Theorem \ref{th : compa rela}, that $(\nu_n)$ is relatively compact, and then converges up to a subsequence (denoted in the same way by $n_k$) to $\phi\in L^{\infty,1}_{mc}([0,T]\times \mathcal{O})^{'}$, for $*\sigma\Big(L^{\infty,1}_{mc}([0,T]\times \mathcal{O})^{'},L^{\infty,1}_{mc}([0,T]\times \mathcal{O})\Big)$, where:
\begin{eqnarray*}
&L^{\infty,1}_{mc}([0,T]\times \mathcal{O}):=\\
&\Big\{~f\in L^\infty([0,T]\times \mathcal{O}),~\exists B\in\mathbb{B}([0,T]),~g\in \mathcal{C}(\mathcal{O}),~f(t,\vartheta)=\mathbf{1}_{B}(t)g(\vartheta)~\Big\}.
\end{eqnarray*}
We now turn to prove that $\phi$ can be represented by an element of $\nu\in \mathcal{M}_{mc}([0,T]\times \mathcal{O})$, that reads: 
\begin{eqnarray*}
\exists\nu\in \mathcal{M}_{mc}([0,T]\times \mathcal{O}),~~\forall f\in L^{\infty,1}_{mc}([0,T]\times \mathcal{O}),~~
\phi(f)=\displaystyle\int_{[0,T]\times \mathcal{O}}f(t,\vartheta)\nu(dt,d\vartheta).
\end{eqnarray*}
For this purpose, we use a Riesz representation Theorem, and more precisely we are going to show that $\phi$ satisfies (i) and (ii) of Theorem \ref{th : Riez}.\\ 
Let $B\in \mathbb{B}([0,T])$, we have:
\begin{eqnarray*}
(t,\vartheta)\mapsto \mathbf{1}_B\otimes 1(t,\vartheta):=\begin{cases}
1,\text{ if }t\in B,\\
0,\text{ if }t\notin B,
\end{cases}
\end{eqnarray*}
belongs to $L^{\infty,1}([0,T]\times \mathcal{O})$ and
\begin{eqnarray*}
\nu_{n_k}(\mathbf{1}_B\otimes 1)\xrightarrow[]{k\to+\infty} \phi(\mathbf{1}_B\otimes 1), \\
\nu_{n_k}(\mathbf{1}_B\otimes 1)=\ell_{n_k}(B)\xrightarrow[]{k\to+\infty} \ell(B).
\end{eqnarray*}
By uniqueness of the weak limit, we get that $\phi(\mathbf{1}_B\otimes 1)=\ell(B)$. Since $\ell\in \mathcal{L}([0,T])$, $\ell$ defines a Borel measure on $([0,T],\mathbb{B}([0,T]))$, implying that (i) of Theorem \ref{th : Riez} holds true. On the other hand, because $\mathcal{O}$ is compact, we deduce that (ii) of Theorem \ref{th : Riez} holds true. It follows that there exists $\nu \in \mathcal{M}_{mc}([0,T]\times \mathcal{O})$, such that:  
\begin{eqnarray*}
\forall f\in L^{\infty,1}_{mc}([0,T]\times \mathcal{O}),~~\phi(f)=\displaystyle\int_{[0,T]\times \mathcal{O}}f(t,\vartheta)\nu(dt,d\vartheta).
\end{eqnarray*}
Since $\phi$ is a continuous linear form on $\text{Span}(L^{\infty,1}_{mc}([0,T]\times \mathcal{O}))$, which is dense in $L^{\infty}_{mc}([0,T]\times \mathcal{O})$ for the uniform convergence (see Lemma \ref{lm : densité }), we deduce that: 
\begin{eqnarray*}
\forall f\in L^\infty_{mc}([0,T]\times \mathcal{O}),~~\phi(f)=\displaystyle\int_{[0,T]\times \mathcal{O}}f(t,\vartheta)\nu(dt,d\vartheta).
\end{eqnarray*}
Finally, to complete the proof, it is enough to show that the projection $\nu^{[0,T]}(dt)$ is equal to $\ell(dt)$. For this we use that for any $B\in\mathbb{B}([0,T])$:
\begin{eqnarray*}
&\displaystyle\int_{[0,T]} \mathbf{1}_{B}(t)\ell_{\nu_{n_k}}(dt)~~\xrightarrow[n_k]{k\to +\infty}~~ \int_{[0,T]}\mathbf{1}_{B}(t)\ell(dt),\\&
\displaystyle\int_{[0,T]\times \mathcal{O}} \mathbf{1}_{B}(t)\nu_{n}(dt,d\vartheta)~~\xrightarrow[n_k]{k\to +\infty}~~\int_{[0,T]\times \mathcal{O}}\mathbf{1}_{B}(t)\nu(dt,d\vartheta).
\end{eqnarray*}
The uniqueness of the weak limit implies:
\begin{eqnarray*}
\forall B\in\mathbb{B}([0,T]),~~\int_{[0,T]}\mathbf{1}_{B}(t)\ell(dt) = \int_{[0,T]\times \mathcal{O}}\mathbf{1}_{B}(t)\nu(dt,d\vartheta)
\end{eqnarray*}
hence:
\begin{eqnarray*}
\ell(dt) = \int_{\mathcal{O}}\nu(dt,d\vartheta),
\end{eqnarray*}
and that completes the proof.
\end{proof}
\begin{Theorem}\label{th : separablitéV} The set
$\mathcal{V}([0,T]\times \mathcal{O})$ of generalized actions at the junction point $\bf 0$, endowed with the weak topology $*\sigma\Big(L^\infty_{mc}([0,T]\times \mathcal{O})^{'},L^\infty_{mc}([0,T]\times \mathcal{O})\Big)$ is Polish. 
\end{Theorem}
\begin{proof}
Recall that $\mathcal{M}_{mc}([0,T]\times \mathcal{O})$ endowed with the weak topology $*\sigma\Big(L^\infty_{mc}([0,T]\times \mathcal{O})^{'},L^\infty_{mc}([0,T]\times \mathcal{O})\Big)$ is separable since:
\begin{eqnarray*}
\mathcal{M}_{mc}([0,T]\times \mathcal{O})~\subset~\bigcup_{n\ge 0}\Big\{~\phi \in L^{\infty}_{mc}([0,T]\times \mathcal{O})^{'},~\|\phi\|\leq n~~\Big\},
\end{eqnarray*}
and from Banach-Alaoglu-Bourbaki's Theorem
\begin{eqnarray*}
\forall n\ge 0,~~\Big\{~\phi \in L^{\infty}_{mc}([0,T]\times \mathcal{O})^{'},~~\|\phi\|\leq n~\Big\}
\end{eqnarray*}
is compact for $*\sigma\Big(L^\infty_{mc}([0,T]\times \mathcal{O})^{'},L^\infty_{mc}([0,T]\times \mathcal{O})\Big)$.\\
As a subset of $\mathcal{M}_{mc}([0,T]\times \mathcal{O})$, we deduce that $\mathcal{V}([0,T]\times \mathcal{O})$ is separable for the weak topology $*\sigma\Big(L^\infty_{mc}([0,T]\times \mathcal{O})^{'},L^\infty_{mc}([0,T]\times \mathcal{O})\Big)$.\\
To conclude, let $(\nu_n(dt,d\vartheta):=\ell_n(dt)\nu_{t,n}(d\vartheta))$ be a Cauchy sequence of $\mathcal{V}([0,T]\times \mathcal{O})$. By definition, we have:
\begin{eqnarray*}
&\forall \varepsilon>0,~~\exists n_0\in \mathbb{N},~~\forall n\ge n_0,~~\forall p\ge 0,~~\forall f\in L^\infty_{mc}([0,T]\times \mathcal{O}),\\&
\Big|~~\displaystyle\int_{[0,T]\times \mathcal{O}}f(t,\vartheta)\nu_{n+p}(dt,d\vartheta)-\displaystyle\int_{[0,T]\times \mathcal{O}}f(t,\vartheta)\nu_{n}(dt,d\vartheta)~~\Big|~~\leq~~\varepsilon.
\end{eqnarray*}
Let $s\in[0,T]$. Choosing $f(t,\vartheta)=\mathbf{1}_{[0,s]}(t)$, we obtain also that $(\ell_n)$ is a Cauchy sequence of $\mathcal{L}([0,T])$ and then converges uniformly to $l\in \mathcal{L}([0,T])$. The converse of Ascoli's Theorem implies that the sequence $(\ell_n)$ satisfies:
\begin{eqnarray*}
&\exists C>0,~~\forall n \ge 0,~~ \ell_n(T)~~\leq~~C,
\\ 
&\exists m\in\mathcal{C}([0,T]),~~m(0)=0,~~\forall n\ge 0,~~\forall (t,s)\in[0,T],~~ |\ell_n(t)-\ell_n(s)|\leq m(|t-s|).
\end{eqnarray*}
Using Theorem \ref{th : compactV}, we deduce that $(\nu_n)$ converges to some $\nu \in \mathcal{V}([0,T]\times \mathcal{O})$ for the weak topology $*\sigma\Big(L^\infty_{mc}([0,T]\times \mathcal{O})^{'},L^\infty_{mc}([0,T]\times \mathcal{O})\Big)$, and that completes the proof.
\end{proof}
\begin{Remark}\label{rm : espace cano Polonais}
As a consequence of Theorem \ref{th : separablitéV}, we obtain that the canonical space:
$$
\Phi~~=~~\mathcal{C}^{\mathcal{N}}([0,T])\times \Big(\prod_{i=1}^I \mathcal{U}([0,T]\times \mathcal{B}_i)\Big)\times \mathcal{V}([0,T]\times \mathcal{O}),
$$ 
involved in the construction of our canonical process $X(\cdot)$ given in \eqref{eq : def processus cano}, is Polish.
\end{Remark}
\subsection{Some estimates and paths properties of the process} \label{subs Some estimations and paths properties of the process}
This subsection is dedicated to give some estimates of the paths for the canonical process $X(\cdot)$. We recall also the crucial “non-stickiness” property at the neighborhood of the junction point $\bf 0$, established in Proposition 4.11 in \cite{Spider}. In order to lighten this contribution,  we do not provide proof of the results that follow. According to assumption $(\mathcal{H})$, the coefficients of diffusion on each ray and the spinning measure at $\bf 0$ are lipschitz continuous uniformly with respect to the control variables. We can then exactly proceed with the same ideas contained in \cite{Spider}.

Adapting the same arguments given in Lemma 4.4 of \cite{Spider}, we can claim:
\begin{Proposition}\label{pr: ineq}
Assume assumption $(\mathcal{H})$. Let $\P^{t,x,i,l}_{\beta,\vartheta}\in \mathcal{A}\big(t,x,i,l\big)$. There exists a uniform constant $C$, depending only on the data $(T,\overline{b},\overline{\sigma})$, such that:
\begin{align*}
\mathbb{E}^{\P^{t,x,i,l}_{\beta,\vartheta}}\Big[~~\big|~x(\cdot)^{2}~\big|_{(t,)}+\big|~\ell_{\nu_0(\cdot)}(\cdot)^{2}~\big|_{(t,)}\Big]~\leq~C(1+x^{2}+l^2),
&\\
\forall \theta \in (0,1),~~\mathbb{E}^{\P^{t,x,i,l}_{\beta,\vartheta}}\Big[~~m\big((x(\cdot),i(\cdot)),\theta\big)^{2}+m\big(\ell_{\nu_0(\cdot)}(\cdot),\theta\big)^{2}\Big]~\leq~C\theta \ln(\frac{2T}{\theta}),
\end{align*}
where recall that: $\forall \big(x(\cdot),i(\cdot)\big)\in \mathcal{C}^{\mathcal{N}}([0,T])$, $\forall \ell(\cdot)\in \mathcal{L}([0,T])$, and $\forall \theta \in (0,1)$:
\begin{align*}
& m\big(\big(x(\cdot),i(\cdot)\big),\theta\big)=\\
&\sup\Big\{~d^{\mathcal{N}}\Big(\big(x(s),i(s)\big),\big(x(u),i(u)\big)\Big),~~(u,s)\in [t,],~~|u-s|\leq \theta,~~\theta\in(0,1)~\Big\},\\ 
&m\big(\ell(\cdot),\theta\big)= \sup \Big\{~|\ell(u)-\ell(s)|,~~(u,s)\in [t,],~~|u-s|\leq \theta,~~\theta\in (0,1)~\Big\}.
\end{align*}
\end{Proposition}
We recall now the central "non-stickiness"  estimate for the process in the neighborhood of the junction point $\bf 0$. In the same spirit of the proof of Proposition 4.11 in \cite{Spider}, we have: 
\begin{Proposition}\label{pr: estimation temps en 0}
Assume assumption $(\mathcal{H})$. Let $\P^{t,x,i,l}_{\beta,\vartheta}\in \mathcal{A}\big(t,x,i,l\big)$. There exists a constant $C>0$, depending only on the data $\Big(T,\overline{b},\overline{\sigma},\underline{\sigma},x\Big)$, introduced in assumption $(\mathcal{H})$, such that:
\begin{eqnarray}\label{eq:majorenzero generale}  \forall \varepsilon>0,~~\mathbb{E}^{\P^{t,x,i,l}_{\beta,\vartheta}}\Big[~~\int_{t}^{T}\mathbf{1}_{\{x(s)<\varepsilon\}}ds~~\Big]~\leq~C\varepsilon.
\end{eqnarray}
\end{Proposition}

We have now the necessary tools to prove the compactness of $\mathcal{A}\big(t,x,i,l\big)$, that is one of the key points in order to show the dynamic programming principle Theorem \ref{th: PPD}.
\begin{Theorem}\label{th: compact}
The set of probability measures $\mathcal{A}\big(t,x,i,l\big)$, endowed with the weak topology, is nonempty and compact.
Moreover, the value function $u$ defined in \eqref{eq : value function} attains its supremum.
\end{Theorem}
\begin{proof}
Recall that the fact that $\mathcal{A}\big(t,x,i,l\big)$ is nonempty is a consequence of the main Theorem in \cite{Spider} (see Remark \ref{rm: non vide}).
Let us show first that $\mathcal{A}\big(t,x,i,l\big)$ is precompact for the weak topology. For this purpose, let $\P_{\beta,\vartheta}^{t,x,i,l} \in \mathcal{A}\big(t,x,i,l\big)$.
It is enough to show that all the following projections: \begin{eqnarray*}
\Big\{~~\P_{\beta,\vartheta}^{t,x,i,l}|_{\mathcal{C}^{\mathcal{N}}([0,T])},~~\P_{\beta,\vartheta}^{t,x,i,l}\in \mathcal{A}\big(t,x,i,l\big)~~\Big\},&\\
\Big(~~\Big\{\P_{\beta,\vartheta}^{t,x,i,l}|_{\mathcal{U}([0,T]\times \mathcal{B}_i)},~~\P_{\beta,\vartheta}^{t,x,i,l}\in \mathcal{A}\big(t,x,i,l\big)~~\Big\}~~\Big)_{i\in\{1,   \ldots    I\}},&\\
\Big\{~~\P_{\beta,\vartheta}^{t,x,i,l}|_{\mathcal{V}([0,T]\times \mathcal{O})},~~\P_{\beta,\vartheta}^{t,x,i,l}\in \mathcal{A}\big(t,x,i,l\big)~~\Big\},&
\end{eqnarray*}
are precompact.

Remark first that the precompactness of $\Big\{~~\P_{\beta,\vartheta}^{t,x,i,l}|_{\mathcal{C}^{\mathcal{N}}([0,T])},~~\P_{\beta,\vartheta}^{t,x,i,l}\in \mathcal{A}\big(t,x,i,l\big)~~\Big\}$ is a consequence of the estimates obtained in Proposition \ref{pr: ineq} and Ascoli's Theorem. On the other hand, since all the sets $(\mathcal{U}([0,T]\times \mathcal{B}_i)_{i\in [I]}$ are compact, we obtain that $\Big(~~\Big\{\P_{\beta,\vartheta}^{t,x,i,l}|_{\mathcal{U}([0,T]\times \mathcal{B}_i)},~~\P_{\beta,\vartheta}^{t,x,i,l}\in \mathcal{A}\big(t,x,i,l\big)~~\Big\}~~\Big)_{i\in [I]}$ are precompact.
We focus in the sequel on the proof of the precompactness for $\Big\{~~\P_{\beta,\vartheta}^{t,x,i,l}|_{\mathcal{V}([0,T]\times \mathcal{O})},~~\P_{\beta,\vartheta}^{t,x,i,l}\in \mathcal{A}\big(t,x,i,l\big)~~\Big\}$.

Let $\varepsilon >0$. It follows from Proposition \ref{pr: ineq}, that there exists a constant $C>0$, depending only on the data $(T,\overline{b},\overline{\sigma})$ such that 
\begin{eqnarray*}
&\mathbb{E}^{\P^{t,x,i,l}_{\beta,\vartheta}}\Big[~~\Big|~~\ell_{\nu_0(\cdot)}(\cdot)^{2}~~\Big|_{(t,)}~~\Big]~~\leq~~C(1+x^2+l^2),\\
&\forall\theta\in(0,T],~~\mathbb{E}^{\P^{t,x,i,l}_{\beta,\vartheta}}\Big[~~m(\ell_{\nu_0(\cdot)}(\cdot),\theta)^{2}~~\Big]~~\leq~~C\theta \ln(\frac{2T}{\theta}).
\end{eqnarray*}
Set: 
\begin{align*}
&K_{\varepsilon}:=\Big\{~~\nu_0\in \mathcal{V}([0,T]\times \mathcal{O}),~~\ell_{\nu_0}(T)\leq \sqrt{\frac{2C(1+x^2+l^2)}{\varepsilon}},
\\ &\forall \theta \in(0,1),~~w(\ell_{\nu_0},\theta)\leq \sqrt{\frac{2C\theta \ln(2T/\theta)}{\varepsilon}}~~\Big\}.
\end{align*}
Using Proposition \ref{th : compactV}, we know that $K_\varepsilon$ is compact for the weak topology $*\sigma\Big(L^\infty_{mc}([0,T]\times \mathcal{O})^{'},L^\infty_{mc}([0,T]\times \mathcal{O})\Big)$. Moreover, using Tchebychev's inequality, we get that
\begin{eqnarray*}
\P_{\beta,\vartheta}^{t,x,i,l}|_{\mathcal{V}([0,T]\times \mathcal{O})}\Big(~~\Big(\nu_0(s)\Big)_{t\leq s\leq T}\notin K_\varepsilon~~\Big)~~\leq~~\varepsilon,
\end{eqnarray*}
which leads to obtain the the precompactness of $\Big\{~~\P_{\beta,\vartheta}^{t,x,i,l}|_{\mathcal{V}([0,T]\times \mathcal{O})},~~\P_{\beta,\vartheta}^{t,x,i,l}\in \mathcal{A}\big(t,x,i,l\big)~~\Big\}$, with the aid of Prokhorov’s theorem and the Polish property established for $\mathcal{V}([0,T]\times \mathcal{O})$ in Theorem \ref{th : separablitéV}. 

Let us show now that $\mathcal{A}\big(t,x,i,l\big)$ is closed for the weak topology. For this purpose, one can consider  $\P_{\beta,\vartheta}^{t,x,i,l}(n) \in \mathcal{A}\big(t,x,i,l\big)$ a sequence converging weakly to $\P_{\beta,\vartheta}^{t,x,i,l}$, and argue exactly with the same arguments given in the proof of Theorem 3.1 (existence part) (see sub Section 5.1 in \cite{Spider}). To conclude, let us introduce the reward function $r_w$, associated to the value function $u$ defined in \eqref{eq : value function}:
\begin{eqnarray*}
&r_w:=\\
&\begin{cases}
\mathcal{A}\big(t,x,i,l\big)\to \R \\
\P^{t,x,i,l}_{\beta,\vartheta}\mapsto \mathbb{E}^{\P^{t,x,i,l}_{\beta,\vartheta}}\Big[~~\displaystyle\displaystyle\int_{t}^{T}\int_{\mathcal{B}_{i(u)}}h_{i(u)}(u,x(u),l(u),\beta_{i(u)})\nu_{i(u)}(T)(du,d\beta_{i(u)})\\
\hspace{2.5cm}+~~\displaystyle \int_{t}^{T}\int_{\mathcal{O}}h_0(u,l(u),\vartheta)\nu_0(T)(du,d\vartheta)~+~g_{i(T)}\big(x(T),\ell(T)\big)~~\Big]
\end{cases}.
\end{eqnarray*}
From assumption $(\mathcal{H})$, the functionals appearing inside of the last expectation are continuous and bounded for the topology induced on the canonical space $\Phi$. Hence, $r_w$ is continuous, and $u$ attains its maximum on the compact set $\mathcal{A}\big(t,x,i,l\big)$.
\end{proof}
\subsection{Dynamic programming principle}\label{sec : DPP}
We conclude this section by the proof of the dynamic programming principle: Theorem \ref{th: PPD}. Most of the main ideas come originally from the works of \cite{ElKaroui}, that use essentially the compactness of the set of admissible rules $\mathcal{A}\big(t,x,i,l\big)$, its stability by conditioning and concatenation, and a measurable lemma selection.

In the sequel, we state two propositions concerning both stability of $\mathcal{A}\big(t,x,i,l\big)$ by conditioning and concatenation.
For the convenience of the reader and in order to make the reading more fluid, we do not sketch the proofs. One can see easily that they are a straight consequence of Theorem 1.3.4, Lemma 6.1.1 and Theorem 6.1.2 in \cite{Stroock}, adapted and valid also for our purpose, because our canonical space $\Phi$ is Polish.
\begin{Proposition}\label{pr: exist traje conca} Let $\tau$ be a $(\Psi_t)_{t\leq s\leq T}$ stopping time, and $(\Q_{\omega})_{\omega \in \Phi}$, a transition probability kernel from $(\Phi,\Psi_{\tau})$ to $(\Phi,\Psi_{T})$ satisfying:
$$\forall \omega\in \Phi,~~\Q_\omega\Big(~~X(\tau(\omega),\cdot)=X(\tau(\omega),\omega)~~\Big)=1.$$
Then:\\
a) there exists a unique transition probability kernel from $(\Phi,\Psi_{\tau})$ to $(\Phi,\Psi_{T})$, denoted by $(\mathbb{I}_\omega\otimes_{\tau(\omega)}\Q_\omega)_{\omega\in \Phi}$, such that:
\begin{eqnarray*}
&\forall \omega\in \Phi,~~\mathbb{I}_\omega\otimes_{\tau(\omega)}\Q_\omega\big(~X\big(s,\cdot\big)=X\big(\tau(\omega),\omega\big),~\forall s\in \big[0,\tau(\omega)\big]~\big)=1,\\&
\forall A\in \sigma\big(~X\big((\tau+s)\wedge T\big),~s\ge 0~\big),~~\mathbb{I}_\omega\otimes_{\tau(\omega)}\Q_\omega(A)=\Q_\omega(A).
\end{eqnarray*}
b) Moreover, if $\P$ is a probability measure on $(\Phi,\Psi_T)$, we associate a unique probability measure on $(\Phi,\Psi_T)$, denoted by $\P\otimes_\tau \Q$ such that:\\
(i) the restriction of $\P\otimes_\tau \Q$ with respect to $\Psi_\tau$ is equal to $\P$,\\
(ii) a r.c.p.d (regular conditional probability distribution) of $\P\otimes_\tau \Q$  with respect to $\Psi_\tau$ is equal to $(\mathbb{I}_\omega\otimes_{\tau(\omega)}\Q_\omega)_{\omega\in \Phi}$.
\end{Proposition}
\begin{Proposition}\label{pr:condi} \textbf{Stability by conditioning:} the set of admissible rules $\mathcal{A}\big(t,x,i,l\big)$ is stable under conditioning, with the following meaning; Let $\P_{\beta,\vartheta}^{t,x,i,l}\in \mathcal{A}\big(t,x,i,l\big)$, and $\tau$ a  $(\Psi_s)_{0\leq s \leq T}$ stopping time, then there exists a probability kernel from $(\Phi,\Psi_{\tau})$ to $(\Phi,\Psi_{T})$ denoted by $\big(\P_{\beta_{\tau(\omega)},\vartheta_{\tau(\omega)}}^{\tau(\omega),x_{\tau(\omega)},i_{\tau(\omega)},l_{\tau(\omega)}}\big)_{\omega\in \Phi}$, satisfying:\\
-there exists $N \subset \Psi_T$, with $\P_{\beta,\vartheta}^{t,x,i,l}(N)=0$ such that:
$$\forall \omega\in \Phi \setminus N,~~\P_{\beta_{\tau(\omega)},\vartheta_{\tau(\omega)}}^{\tau(\omega),x_{\tau(\omega)},i_{\tau(\omega)},l_{\tau(\omega)}} \in \mathcal{A}\Big(\tau(\omega),x_{\tau(\omega)},i_{\tau (\omega)},l_{\tau(\omega)}\Big),$$
-for all $f:\Phi\to \R$, and $\sigma\big(~X\big((\tau+s)\wedge T\big),~s\ge 0~\big)$ measurable, we have: $$\displaystyle \E^{\P_{\beta,\vartheta}^{t,x,i,l}}\big[~f~\big|~\Psi_\tau~\big]=\displaystyle\E^{ \P_{\beta_{\tau(\omega)},\vartheta_{\tau(\omega)}}^{\tau(\omega),x_{\tau(\omega)},i_{\tau(\omega)},l_{\tau(\omega)}}}\big[~f~\big],~~\P_{\beta,\vartheta}^{t,x,i,l}~~\text{a.s.}$$ 
\end{Proposition}
\begin{Proposition}\label{pr:conca} \textbf{Stability by concatenation:} the set of admissible rules
$\mathcal{A}\big(t,x,i,l\big)$ is stable under concatenation, with the following meaning: Let $\P_{\beta,\vartheta}^{t,x,i,l}\in \mathcal{A}\big(t,x,i,l\big)$ and $\tau$ a $(\Psi_s)_{t\leq s\leq T}$ stopping time. Assume that $\big(\P_{\beta_{\tau(\omega)},\vartheta_{\tau(\omega)}}^{\tau(\omega),x_{\tau(\omega)},i_{\tau(\omega)},l_{\tau(\omega)}}\big)_{\omega \in \Phi}$ is a transition probability kernel from $(\Omega,\Psi_{\tau})$ to $(\Omega,\Psi_{T})$ satisfying:
$$\P_{\beta_{\tau(\omega)},\vartheta_{\tau(\omega)}}^{\tau(\omega),x_{\tau(\omega)},i_{\tau(\omega)},l_{\tau(\omega)}} \in \mathcal{A}\Big(\tau(\omega),x_{\tau(\omega)},i_{\tau (\omega)},l_{\tau(\omega)}\Big),~~\P_{\beta,\vartheta}^{t,x,i,l}~~\text{a.s.}$$
Then $\P_{\beta,\vartheta}^{t,x,i,l}\otimes_\tau \P_{\beta_\tau,\vartheta_\tau}^{\tau,x_\tau,i_\tau,l_\tau} \in \mathcal{A}\big(t,x,i,l\big)$, where $\P_{\beta,\vartheta}^{t,x,i,l}\otimes_\tau \P_{\beta_\tau,\vartheta_\tau}^{\tau,x_\tau,i_\tau,l_\tau}$ is defined in Proposition \ref{pr: exist traje conca} b).
\end{Proposition}
\subsection{Proof of Theorem \ref{th: PPD}} 
\begin{proof}
One can exactly follows the same lines given in the proof of Theorem 6.3 (point c) given in \cite{ElKaroui}, since the set of admissible rules $\mathcal{A}\big(t,x,i,l\big)$ is:\\
-stable by conditioning (Proposition \ref{pr:conca}),\\
-stable by concatenation (Proposition \ref{pr:conca}),\\
whereas the canonical space $\Phi$ is Polish. Note that also the adaptation of the  measurable selection Lemma (recalled in Proposition \ref{cr: sele mesur}) still holds true for our purpose.
\end{proof}
\section{Comparison theorem for continuous viscosity backward parabolic Walsh's spider HJB equations having a non linear local-time Kirchhoff's boundary condition}\label{sec : theo compa viscosity}
As announced in Introduction, the first result giving a comparison theorem for fully non linear HJB equations of second order, with non vanishing viscosity at the vertex, and having a non linear Kirchhoff's boundary condition at $\bf 0$ was proved in \cite{Ohavi Walsh PDE}. More precisely, therein the system called: {\it Walsh’s spider Hamilton-Jacobi-Bellman}, involves a more sophisticated new boundary condition named: {\it non linear local-time Kirchhoff’s transmission at $\bf 0$}, (in the elliptic case and  for a bounded domain). Let us describe the new ingredients that were used in \cite{Ohavi Walsh PDE} to obtain the central comparison theorem. First note that it is easy to be convinced that at the interior of each edge, one can argue with classical arguments (doubling variable method), as soon as the Hamiltonians satisfy the classical assumptions. The main central difficulties are obviously concentrated at $\bf 0$ because of the discontinuities of the Hamiltonians, and the presence of the Kirchhoff's boundary condition. The introduction of an external deterministic ’local-time’ variable $l$ in \cite{Ohavi Walsh PDE} – that is the counterpart to the local time $\ell$ – is one of the main crucial
ingredient to obtain the comparison principle. More precisely, the main idea was to build test functions at the neighborhood of the vertex solutions of ODE, with well-designed coefficients depending on what the author calls {\it the speed of the Hamiltonian} (by analogy to the speed measure of a diffusion). The crucial key point therein was to impose a 'local-time' derivative at the vertex absorbing the error term induced by the {\it Kirchhoff's speed of the Hamiltonians}; that is the Kirchhoff's average of all the speed of the Hamiltonians at the neighborhood of $\bf 0$. Note that even without the presence of
the external variable $l$ in the original HJB problem, the ’artificial’ introduction of this external variable in the problem allowed to extend the main results contained in \cite{Lions Souganidis 1} and \cite{Lions Souganidis 2} to the fully non linear and non degenerate framework. Finally, it is important also to emphasize that the comparison theorem for viscosity solutions in \cite{Ohavi Walsh PDE} is stated in the strong sense for the Kirchhoff's boundary condition at $\bf 0$, namely without any dependency of the values of the Hamiltonians at $\bf 0$, which is also innovative for Neumann problems in the non linear case.

To be as clear as possible, since the ideas introduced in \cite{Ohavi Walsh PDE} are quite new in the field of viscosity theory, in what follows we give the main lines of the proof for Theorem \ref{th compa visco}, since the central ideas are contained in Theorem 2.2 of \cite{Ohavi Walsh PDE}.

\textbf{Proof of Theorem \ref{th compa visco}.} 
\begin{proof}
Fix $\lambda>0$. Let $\mathfrak{f}$ be a super solution of \eqref{eq PDE Walsh} and  $\mathfrak{g}$ be a sub solution of \eqref{eq PDE Walsh}, both uniformly bounded and satisfying the following backward boundary condition $\mathfrak{f}_{\cdot}(T,\cdot,\cdot)\ge \mathfrak{g}_{\cdot}(T,\cdot,\cdot)$. Using the Definition \ref{def viscosity}, it is easy to obtain that the following map:
$$f_i(t,x,l):=f^\lambda_i(t,x,l)=e^{\lambda t}\mathfrak{f}_i(t,x,l),~~(i,t,x,l)\in [I]\times [0,T]\times[0,+\infty)^2,$$
(resp. $g_i(t,x,l):=g^\lambda_i(t,x,l)=e^{\lambda t}\mathfrak{g}_i(t,x,l)$) is super solution (resp. sub solution) of the following HJB system posed on the domain $\mathcal{D}_T$:
\begin{eqnarray}\label{eq PDE Walsh 2}
\begin{cases}
\textbf{HJB equation parameterized by the local-time on each ray :}\\
\partial_t \mathfrak{p}_i(t,x,l)-\lambda \mathfrak{p}_i(t,x,l)+\underset{\beta_i\in \mathcal{B}_i}{\sup}\Big\{\ds \frac{1}{2}\sigma_i^2(t,x,l,\beta_i)\partial^2_x\mathfrak{p}_i(t,x,l)+\\
b_i(t,x,l,\beta_i)\partial_x\mathfrak{p}_i(t,x,l)+e^{\lambda t}h_i(t,x,l,\beta_i)\Big\}=0,~~(t,x,l)\in(0,T)\times (0,+\infty)^2,\\
\textbf{Non linear local-time Kirchhoff's boundary transmission :}\\
\partial_l\mathfrak{p}(t,0,l)+\underset{ \vartheta \in \mathcal{O}}{\sup} \Big\{\ds\sum_{i=1}^I\mathcal{S}_i(t,l,\vartheta)\partial_x\mathfrak{p}_i(t,0,l)+\\
\hspace{3 cm}e^{\lambda t}h_0(t,l,\vartheta)\Big\}=0,~~(t,l)\in(0,T)\times(0,+\infty)\\
\textbf{Continuity condition at $\bf 0$} :\\
\forall (i,j)\in[I]^2,~~\forall (t,l)\in [0,+\infty),~~\mathfrak{p}_i(t,0,l)=\mathfrak{p}_j(t,0,l).
\end{cases}
\end{eqnarray}
in the sense of Definition \ref{def viscosity}, satisfying the backward boundary condition: $f_{\cdot}(T,\cdot,\cdot)\ge g_{\cdot}(T,\cdot,\cdot)$.
It appears clear therefore to obtain a comparison Theorem for the HJB system \eqref{th compa visco}, it is enough to obtain a comparison theorem for the modified HJB system \eqref{eq PDE Walsh 2}. Let then $f$ be a super solution (resp. $g$ a sub solution) of the last system \eqref{eq PDE Walsh 2} satisfying $f_{\cdot}(T,\cdot,\cdot)\ge g_{\cdot}(T,\cdot,\cdot)$.\\
Fix $\alpha>0$ and $(\underline{s},\underline{\ell}) \in (0,T)\times (0,+\infty)$.
We argue by contradiction assuming that:
\begin{eqnarray}\label{eq: contra compa}
 \nonumber &\sup\Big\{~~ g_{i}(t,x,l)- f_{i}(t,x,l)-\alpha(x^2+l^2),\\
 &(i,t,x,l)\in [I]\times [\underline{s},T]\times[0,+\infty)\times[\underline{\ell},+\infty)~~\Big\}~>0.
\end{eqnarray}
Since $g$ and $f$ are uniformly bounded and continuous, the last supremum is reached at some point $(i_\star,t_\star,x_\star,l_\star)\in [I]\times [\underline{s},T)\times[0,+\infty)\times[\underline{\ell},+\infty)$ such that:
\begin{eqnarray}\label{eq positive contra sur sous solu}
 g_{i_\star}(t_\star,x_\star,l_\star)- f_{i_\star}(t_\star,x_\star,l_\star)>0. 
\end{eqnarray}
\textbf{Step 1:} Assume that $x_\star>0$. Denote $H_{i_\star}^{l_\star}$ the following Hamiltonian acting on the ray $\mathcal{R}_{i_\star}$:
$$H_{i_\star}^{l_\star}:=\begin{cases}
[\underline{s},T]\times\R\times[0+\infty)\times \R^3\to \R\\
 (s,q,x,u,p.S)\mapsto q-\lambda u+\underset{\beta_{i_\star} \in \mathcal{B}_{i_\star}}{\sup}\Big\{\ds\frac{1}{2}\sigma_{i_\star}^2(s,x,l_\star,\beta_{i_\star})S~+\\
 \hspace{3.8cm}b_{i_\star}(s,x,l_\star,\beta_{i_\star})p+e^{\lambda t}h_{i_\star}(s,x,l_\star,\beta_{i_\star})\Big\}  
\end{cases},$$
where the deterministic level of the local-time variable is fixed at $l=l_\star$. As explained in the proof of Theorem 2.2 (Step 1), one can proceed using classical arguments, because the last Hamiltonian $H_{i_\star}^{l_\star}$ does not have any dependency with some derivative w.r.t the variable $l$. Let us briefly recall the main steps, coming from the seminal survey on viscosity theory solutions, (Theorem 8.2 in \cite{User guide}):\\
-one can use the doubling variable on a thin neighborhood $[t_\star-\eta,t_\star+\eta]\times\mathcal{V}_{i_\star}(x_\star)\ni (u,x),~\eta>0$, such that our space variable $x\in \mathcal{V}_{i_\star}(x_\star)$ belongs included at the interior of the ray $\mathcal{R}_{i_\star}$. More precisely, one has to introduce the following map (here parameterized by $(i_\star.l_\star)$), defined $\forall \varepsilon>0$ by:
\begin{eqnarray*}
 &\mathfrak{w}_{\varepsilon,\eta}^{(i_\star.l_\star)}(t,x,y)=g_{i_\star}(t,x,l_\star)-f_{i_\star}(t,y,l_\star)-\ds \frac{1}{2\varepsilon^2}|x-y|^2,\\
 &(t,x,y)\in [t_\star-\eta,t_\star+\eta]\times\mathcal{V}_{i_\star}(x_\star)^2.   \end{eqnarray*}
-obtain a contradiction with the aid of the Ishii's matrix lemma for the parabolic framework (see for example Theorem 8.3 in \cite{User guide}). Recall that  from assumption $(\mathcal{H})$, it follows that $H_{i_{\star}}^{l_\star}$ is continuous, satisfies the structure growth given in Example 3.6 of \cite{User guide}, whereas the equivalent definition of a super (resp. sub solutions) with the closure of the second-order super jet (resp. subjet) for $f$ (resp. $g$) at $(t_\star,x_\star)$ holds true.

\textbf{Step 2 :} Assume that $x_\star=0$. Here we will give the main ideas and steps that will lead to a contradiction with \eqref{eq positive contra sur sous solu}. We will try as much as possible not to overburden the calculations to make the reading clear. Indeed, the only difference with the proof of the central Theorem 2.2 (Steps 2-3-4) given in \cite{Ohavi Walsh PDE}, is that this time we have the presence of the time variable, but this will have no impact, because the non linear Kirchhoff's condition does not present any dependency with some derivative w.r.t to the time variable. 

In the sequel, we argue like in Theorem 2.2 (Step 2) of \cite{Ohavi Walsh PDE} adapted to our framework. 
We scale first $f$ and $g$ at the vertex $\bf0$, setting 
$$\Theta(f,g)=\frac{1}{2}\big(f(t_\star,0,l_\star)+g(t_\star,0,l_\star)\big).$$
Therefore:
\begin{eqnarray}\label{eq fonction scaller}
 u=f-\Theta(f,g),~~v=g-\Theta(f,g),  
\end{eqnarray}
are respectively super solution and sub solution of the following system with {\it non linear local time's Kirchhoff's boundary transmission}, posed on the domain $\mathcal{D}_T$:
\begin{eqnarray}\label{eq PDE Walsh 2 scaller}
\begin{cases}
\textbf{HJB equation parameterized by the local-time on each ray :}\\
\partial_t \mathfrak{p}_i(t,x,l)-\lambda \mathfrak{p}_i(t,x,l)-\lambda \Theta(f,g)+\underset{\beta_i\in \mathcal{B}_i}{\sup}\Big\{\ds \frac{1}{2}\sigma_i^2(t,x,l,\beta_i)\partial^2_x \mathfrak{p}_i(t,x,l)+\\
b_i(t,x,l,\beta_i)\partial_x\mathfrak{p}_i(t,x,l)+e^{\lambda t}h_i(t,x,l,\beta_i)\Big\}=0,~~(t,x,l)\in(0,T)\times (0,+\infty)^2,\\
\textbf{Non linear local-time Kirchhoff's boundary transmission :}\\
\partial_l\mathfrak{p}(t,0,l)+
\underset{ \vartheta \in \mathcal{O}}{\sup} \Big\{\ds\sum_{i=1}^I\mathcal{S}_i(t,l,\vartheta)\partial_x\mathfrak{p}_i(t,0,l)+\\
\hspace{3 cm}e^{\lambda t}h_0(t,l,\vartheta)\Big\}=0,~~(t,l)\in(0,T)\times(0,+\infty)\\
\textbf{Continuity condition at $\bf 0$} :\\
\forall (i,j)\in[I]^2,~~\forall (t,l)\in [0,+\infty),~~\mathfrak{p}_i(t,0,l)=\mathfrak{p}_j(t,0,l),
\end{cases},~~
\end{eqnarray}
satisfying the following backward boundary condition $u_{\cdot}(T,\cdot,\cdot)\ge v_{\cdot}(T,\cdot,\cdot)$.

Remark this important facts:
\begin{eqnarray}\label{eq positiv u v}
\nonumber &u(t_\star,0,l_\star)=\frac{1}{2}\big(f(t_\star,0,l_\star)-g(t_\star,0,l_\star)\big)<0,\\
\nonumber &v(t_\star,0,l_\star)=\frac{1}{2}\big(g(t_\star,0,l_\star)-f(t_\star,0,l_\star)\big)>0,\\
&v(t_\star,0,l_\star)-u(t_\star,0,l_\star)=g(t_\star,0,l_\star)-f(t_\star,0,l_\star)>0.
\end{eqnarray}
Using \eqref{eq positiv u v} and the continuity of $u$ and $v$, there exists a neighborhood of the vertex $(t_\star,0,l_\star)$ with respect to the variables $(t,x,l)$, denoted by $\mathcal{V}\Big((t_\star,0,l_\star),\big(\varepsilon,\kappa\big)\Big)$ such that (where $\varepsilon>0$, $\kappa>0$ are two small enough parameters expected to be sent to $0$):
\begin{eqnarray*}
&\mathcal{V}\Big((t_\star,0,l_\star),\big(\varepsilon,\kappa\big)\Big):=\Big\{\big(t,(x,i),l\big)\in (0,T)\times\mathcal{N}\times (0,+\infty),~~0\leq x\leq \varepsilon,\\
&0<t_\star-\kappa\leq t \leq  t_\star+\kappa<T,~~0<l_\star-\kappa\leq l \leq  l_\star+\kappa\Big\},   
\end{eqnarray*}
and:
\begin{eqnarray}\label{eq positi u et v}
&\nonumber \forall\big(t,(x,i),l\big)\in \mathcal{V}\Big((t,0,l_\star),\big(\varepsilon,\kappa\big)\Big):~~v_i(t,x,l)\ge 0,~~u_i(t,x,l)\leq 0,\\
&v_i(t,x,l)-u_i(t,x,l)\ge 0.
\end{eqnarray}
In the sequel we introduce also $\eta>0$ and $\gamma>0$ two small strictly positive parameters, designed to drive the construction of the test functions at the neighborhood of $(t_\star,0,l_\star)$. Like in the Step 2 of Theorem 2.2 in \cite{Ohavi Walsh PDE} we introduce the following quantities (recall that $\underline{\sigma}$ states for the ellipticity constant of the second order terms, and $\overline{\sigma}$ their uniform upper bound - see assumption $(\mathcal{H})$);
\begin{eqnarray}\label{eq  expression const scal}
\nonumber \overline{\rho}\big(\lambda,\Theta(f,g)\big)=\lambda\Theta(f,g)\big(\frac{1}{\underline{\sigma}^2}\mathbf{1}_{\Theta(f,g)>0}+\frac{1}{\overline{\sigma}^2}\mathbf{1}_{\Theta(f,g)\leq 0}\big),\\
\nonumber \underline{\rho}\big(\lambda,\Theta(f,g)\big)=\lambda\Theta(f,g)\big(\frac{1}{\underline{\sigma}^2}\mathbf{1}_{\Theta(f,g)\leq 0}+\frac{1}{\overline{\sigma}^2}\mathbf{1}_{\Theta(f,g)> 0}\big),\\
\nonumber \overline{u}^\kappa(0)=\sup\big\{u(t,0,l),~~(t,l)\in[t_\star-\kappa,t_\star+\kappa]\times[l_\star-\kappa,l_\star+\kappa]\big\},\\
\nonumber\forall i\in [I],~~
\underline{u}_i^\kappa(\varepsilon)=\inf\big\{u_i(t,\varepsilon,l),~~(t,l)\in[t_\star-\kappa,t_\star+\kappa]\times[l_\star-\kappa,l_\star+\kappa]\big\},\\
\nonumber \underline{v}^\kappa(0)=\inf\big\{v(t,0,l),~~(t,l)\in[t_\star-\kappa,t_\star+\kappa]\times[l_\star-\kappa,l_\star+\kappa]\big\},\\  \forall i\in [I],~~
\overline{v}_i^\kappa(\varepsilon)=\sup\big\{v_i(t,\varepsilon,l),~~(t,l)\in[t_\star-\kappa,t_\star+\kappa]\times[l_\star-\kappa,l_\star+\kappa]\big\},
\end{eqnarray}
Let $\overline{S}\ge 0$ and $\underline{S}\ge 0$ be two parameters designed to characterize the 'local-time' derivative of the test functions at $\bf 0$.
Using Proposition 4.1 given in \cite{Ohavi Walsh PDE}, we introduce \begin{align*}
&\overline{\phi}=\overline{\phi}(u,\varepsilon,\eta,\kappa,\gamma,\overline{S}),~~\underline{\phi}=\underline{\phi}(v,\varepsilon,\kappa,\eta,\gamma,\underline{S})
\end{align*}
(denoted in the next lines $\big(\overline{\phi},\underline{\phi}\big)$ for the seek of clarity) the two solutions of the  two following ordinary parametric differential equation systems posed on $\mathcal{N}_\varepsilon\times [l_\star-\kappa,l_\star+\kappa]$, independent of the time variable, and parameterized by the local time variable $l$:
\begin{eqnarray}\label{expression EDO sur sol}
\begin{cases}
\ds \frac{2\lambda }{\overline{\sigma}^2}\overline{\phi}_i-\partial_x^2\overline{\phi}_i(x,l)+2\overline{\rho}\big(\lambda,\Theta(f,g)\big)\\
 \ds \hspace{1,0 cm}+2\Big(|b||\partial_x\overline{\phi}_i(x,l)|+e^{\lambda T}|h|\Big)\big/\underline{\sigma}^2=-\eta,~~(x,l)\in (0,\varepsilon)\times(l_\star-\kappa,l_\star+\kappa),\\
\overline{\phi}(0,l)=\overline{u}^\kappa(0)+\overline{S}(l-l_\star),\\
\overline{\phi}_i(\varepsilon,l)=\underline{u}_i^\kappa(\varepsilon)-\gamma+\overline{S}(l-l_\star),~~l\in[l_\star-\kappa,l_\star+\kappa],~~i\in [I].
\end{cases},
\end{eqnarray}
and:
\begin{eqnarray} \label{expression EDO sous sol}
\begin{cases}
\ds \frac{2\lambda }{\overline{\sigma^2}}\underline{\phi}_i-\partial_x^2\underline{\phi}_i(x,l)+2\underline{\rho}\big(\lambda,\Theta(f,g)\big)\\
 \ds \hspace{1,0 cm}-2\Big(|b||\partial_x\underline{\phi}_i(x,l)|+e^{\lambda T}|h|\Big)\big/\underline{\sigma^2}=\eta,~~(x,l)\in (0,\varepsilon)\times(l_\star-\kappa,l_\star+\kappa),\\
\underline{\phi}(0,l)=\underline{v}^\kappa(0)-\underline{S}(l-l_\star),\\
\underline{\phi}_i(\varepsilon,l)=\overline{v}_i^\kappa(\varepsilon)+\gamma-\underline{S}(l-l_\star),~~l\in[l_\star-\kappa,l_\star+\kappa],~~i\in [I].
\end{cases}.
\end{eqnarray}
Assumption $(\mathcal{H})$ combined with Proposition 4.1 given \cite{Ohavi Walsh PDE}, state that both solutions $\overline{\phi}$ and $\underline{\phi}$ are unique and in the class of test functions $\mathcal{C}^{1,2,0}_{{\bf  0},1}\big(\mathcal{N}_\varepsilon\times [l_\star-\kappa,l_\star+\kappa]\big)$, and therefore are test viscosity functions of system \eqref{eq PDE Walsh 2 scaller}, in the sense of Definition \ref{def viscosity}.

Next with the same arguments given in Step 3 Theorem 2.2 of \cite{Ohavi Walsh PDE}, it is easy to check that
$\overline{\phi}$ (resp. $\underline{\phi}$) is a test function of the super solution $u$ (resp. the sub solution $v$) at the vertex $\bf 0$ of the Walsh's spider HJB system \eqref{eq PDE Walsh 2 scaller}.
Therefore there exists $(\overline{t}_\kappa,\overline{l}_\kappa)\in [t_\star-\kappa,t_\star+\kappa]\times[l_\star-\kappa,l_\star+\kappa]$ such that:
\begin{eqnarray}\label{eq second jet u}
\partial_l\overline{\phi}(0,\overline{l}_\kappa)+\sup_{\vartheta \in \mathcal{O}}\Big\{\ds \sum_{i=1}^I\mathcal{S}_i(\overline{t}_\kappa,\overline{l}_\kappa,\vartheta)\partial_x\overline{\phi}_i(0,\overline{l}_\kappa)+e^{\lambda \overline{t}_\kappa}h_0(\overline{t}_\kappa,\overline{l}_\kappa,\vartheta)\Big\}\ge 0, 
\end{eqnarray}
and $(\underline{t}_\kappa,\underline{l}_\kappa)\in [t_\star-\kappa,t_\star+\kappa]\times[l_\star-\kappa,l_\star+\kappa]$ satisfying:
\begin{eqnarray}\label{eq second jet v}
\partial_l\underline{\phi}(0,\underline{l}_\kappa)+\sup_{\vartheta \in \mathcal{O}}\Big\{\ds \sum_{i=1}^I\mathcal{S}_i(\underline{t}_\kappa,\underline{l}_\kappa,\vartheta)\partial_x\underline{\phi}_i(0,\underline{l}_\kappa)+e^{\lambda \underline{t}_\kappa}h_0(\underline{l}_\kappa,\vartheta)\Big\}\leq 0.  
\end{eqnarray}
To conclude, we follow the ideas and calculations of Step 4 of Theorem 2.2 given in \cite{Ohavi Walsh PDE}. Without loss of generality, we know from Proposition 4.1 in \cite{Ohavi Walsh PDE}, that for $\varepsilon<<1,~\kappa<<1$, there exists 
\begin{align*}
&\overline{S}=\overline{S}(\overline{\zeta},\varepsilon,\kappa,\eta,\gamma)=\overline{S}\Big(\overline{\zeta},\varepsilon,\kappa,|\overline{\rho}\big(\lambda, \Theta(f,g)\big)|+\eta,|h|,|\partial_x\overline{\phi}|,\overline{u}^\kappa(0),\big(\underline{u}_i^\kappa(\varepsilon)\big)_{i\in [I]}\Big),~~\\
&\underline{S}=\underline{S}(\overline{\zeta},\varepsilon,\kappa,\eta,\gamma)=\overline{S}\Big(\overline{\zeta},\varepsilon,\kappa,|\underline{\rho}\big(\lambda, \Theta(f,g)\big)|+\eta,|h|,|\partial_x\underline{\phi}|,\underline{v}^\kappa(0),\big(\overline{v}_i^\kappa(\varepsilon)\big)_{i\in [I]}\Big),  
\end{align*}
such that the 'local-time' derivatives of the test functions $\big(\overline{\phi},\underline{\phi}\big)$ satisfy:
\begin{eqnarray*}
&\label{eq absor derivee 1}
\partial_l\overline{\phi}(0,\overline{l}_\kappa)=\overline{S}(\overline{\zeta},\varepsilon,\kappa,\eta,\gamma)\ge \varepsilon I\overline{\zeta}\Big(|\overline{\rho}\big(\lambda, \Theta(f,g)\big)|+\eta+\frac{|b||\partial_x\overline{\phi}|+e^{\lambda T}|h|}{\underline{\sigma}}\Big),\\ &\label{eq absor derivee 2}
\partial_l\underline{\phi}(0,\underline{l}_\kappa)=-\underline{S}(\overline{\zeta},\varepsilon,\kappa,\eta,\gamma)\leq -\varepsilon I\overline{\zeta}\Big(|\underline{\rho}\big(\lambda, \Theta(f,g)\big)|+\eta+\frac{|b||\partial_x\underline{\phi}|+e^{\lambda T}|h|}{\underline{\sigma}}\Big).
\end{eqnarray*}
Recall that $\overline{\zeta}$ is given in assumption $(\mathcal{H}-\bf R)$ and we have that:
\begin{eqnarray}\label{eq cv point}
 \lim_{\kappa \searrow 0} \overline{l}_\kappa =l_\star,~~ \lim_{\kappa \searrow 0} \underline{l}_\kappa =l_\star,~~\lim_{\kappa \searrow 0} \overline{t}_\kappa =t_\star,~~ \lim_{\kappa \searrow 0} \underline{t}_\kappa =t_\star,
\end{eqnarray}
whereas from the expressions given in \eqref{eq  expression const scal} and the continuity of $u$ and $v$:
\begin{eqnarray}\label{eq cv point bord}
&\nonumber \forall i\in[I],~~\lim_{\varepsilon \searrow 0}\limsup_{\kappa \searrow 0} \overline{v}^\kappa_i(\varepsilon) =v(t_\star0,l_\star),~~\lim_{\varepsilon \searrow 0}\limsup_{\kappa \searrow 0} \underline{u}^\kappa_i(\varepsilon) =u(t_\star,0,l_\star),\\
&\lim_{\kappa \searrow 0} \underline{v}^\kappa(0) =v(t_\star,0,l_\star),~~\lim_{\kappa \searrow 0} \overline{u}^\kappa(0) =u(t_\star,0,l_\star).
\end{eqnarray}
To conclude, one can follow the same ideas and calculations in Step 4 of \cite{Ohavi Walsh PDE} to get:
$$0\ge g(t_\star,0,l_\star)-f(t_\star,0,l_\star),$$
which leads to a contradiction with \eqref{eq: contra compa}. Hence, for all $\underline{s}\in (0,T)$, for all $\underline{\ell}\in (0,+\infty)$, for all $t\in [\underline{s},T]$, for all $(x,i)\in \mathcal{N}$ and for all $l\in [\underline{\ell},+\infty)$, we obtain that (after sending $\alpha \searrow 0$ in \eqref{eq: contra compa}):
$$f_i(t,x,l)\ge g_i(t,x,l).$$
We can then complete the proof using the continuity of $f$ and $g$ with respect to the variables $t$ and $l$.
\end{proof}
\section{Unique characterization of the value function for the problem of stochastic scattering control}\label{sec: caracterisation value function}
In this Section, we show that the value function defined in Definition \ref{def fonction valeur}-\eqref{eq : value function}, is characterized in a unique way with the aid of HJB system \eqref{eq PDE Walsh}.
\begin{Proposition}\label{pr: continuite fonction valeur}
Assume assumption $(\mathcal{H})$. Then the value function $u$ defined in Definition \ref{def fonction valeur} - \eqref{eq : value function} is continuous.  
\end{Proposition}
\begin{proof}
 Let $(t,x,i,l)\in [0,T]\times[0,+\infty)\times[I]\times [0,+\infty)$ and $(t_n,x_n,i_n,l_n)_{n\ge 0}$ a sequence converging to $(t,x,i,l)$ for the usual metric defined on $[0,T]\times\mathcal{N}\times[0,+\infty)$ by: $|\cdot|+d^{\mathcal{J}}+|\cdot|$. 
First using assumption $(\mathcal{H})$ and especially the fact that the coefficients of diffusion on each ray and the spinning measure at $\bf 0$ are uniformly Lipschitz continuous w.r.t to the control variables, which ensures the uniqueness in distribution of any $\P^{t,x,i,l}_{\beta,\vartheta}\in \mathcal{A}\big(t,x,i,l\big)$; one can exactly argue like in the proof of Proposition 3.4 of \cite{Spider} to show that the following map:
$$
\begin{cases}
[0,T]\times[0,+\infty)\times[I]\times [0,+\infty)\to \mathcal{P}\big(\Phi,\Psi\big),\\ 
(t,x,i,l) \to \P^{t,x,i,l}_{\beta,\vartheta}
\end{cases}
$$
is continuous. It follows therefore that the value function $u$ defined in Definition \ref{def fonction valeur} - \eqref{eq : value function} is lower semi continuous.

On the other hand, from Theorem \ref{th: compact}, we know that the value function $(t_n,x_n,i_n,l_n)\mapsto u_{i_n}(t_n,x_n,l_n)$ attains its supremum for a given $\overline{\P}^{t_n,x_n,i_n,l_n}_{\beta_n,\mathcal{O}_n}$ belonging to the set of admissible rues $\mathcal{A}\big(t_n,x_n,i_n,l_n\big)$, such that:
 \begin{eqnarray}\label{eq lim u}
\nonumber u_{i_n}(t_n,x_n,l_n)=\mathbb{E}^{\overline{\P}^{t_n,x_n,i_n,l_n}_{\beta_n,\mathcal{O}_n}}\Big[~~\displaystyle\displaystyle\int_{t}^{T}\int_{\mathcal{B}_{i(u)}}h_{i(u)}(u,x(u),l(u),\beta_{i(u)})\nu_{i(u)}(T)(du,d\beta_{i(u)})+\\
\displaystyle \int_{t}^{T}\int_{\mathcal{O}}h_0(u,l(u),\vartheta)\nu_0(T)(du,d\vartheta)+g_{i(T)}\big(x(T),\ell(T)\big)~~\Big].
\end{eqnarray}
Arguing like in the proof of Theorem \ref{th: compact}, we get that the sequence $(\overline{\P}^{t_n,x_n,i_n,l_n}_{\beta_n,\mathcal{O}_n})_{n\ge 0}$ is tight, and therefore converges weakly up to a bus sequence $n_k$ to some probability measure $\overline{\P}\in \mathcal{A}\big(t,x,i,l\big)$. Once again, with the same arguments of Proposition 3.4 in \cite{Spider}, we have that $\overline{\P}$ is solution of the spider martingale problem $\big(\mathcal{S}_{pi}-\mathcal{M}_{ar}\big)$, satisfying the conditions $(\mathcal{S})$, associated with the initial condition $(t,x,i,l)$. Denote then in the sequel $\overline{\P}=\overline{\P}^{t,x,i,l}_{\beta,\mathcal{O}}$.
Sending $k\to +\infty$ in \eqref{eq lim u}, it follows that:
 \begin{align*}
\lim_{n_k\to +\infty}u_{i_{n_k}}(t_{n_k},x_{n_k},l_{n_k})=\mathbb{E}^{\overline{\P}^{t,x,i,l}_{\beta,\mathcal{O}}}\Big[~~\displaystyle\displaystyle\int_{t}^{T}\int_{\mathcal{B}_{i(u)}}h_{i(u)}(u,x(u),l(u),\beta_{i(u)})\nu_{i(u)}(T)(du,d\beta_{i(u)})+\\
\displaystyle \int_{t}^{T}\int_{\mathcal{O}}h_0(u,l(u),\vartheta)\nu_0(T)(du,d\vartheta)+g_{i(T)}\big(x(T),\ell(T)\big)~~\Big].
\end{align*}
Using one more time the uniqueness in distribution for the martingale problem $\big(\mathcal{S}_{pi}-\mathcal{M}_{ar}\big)$, satisfying the conditions $(\mathcal{S})$, and the fact that the weak convergence of probability measure is metrizable, it follows that $\overline{\P}^{t,x,i,l}_{\beta,\mathcal{O}}$ is the unique limit point of $(\overline{\P}^{t_n,x_n,i_n,l_n}_{\beta_n,\mathcal{O}_n})_{n\ge 0}$. We conclude that $(\overline{\P}^{t_n,x_n,i_n,l_n}_{\beta_n,\mathcal{O}_n})_{n\ge 0}$ converges to $\overline{\P}^{t,x,i,l}_{\beta,\mathcal{O}}$ and then:
 \begin{eqnarray*}
&\ds \underset{n\to +\infty}{\overline{\lim}}u_{i_{n}}(t_{n},x_{n},l_{n})=\mathbb{E}^{\overline{\P}^{t,x,i,l}_{\beta,\mathcal{O}}}\Big[~~\displaystyle\displaystyle\int_{t}^{T}\int_{\mathcal{B}_{i(u)}}h_{i(u)}(u,x(u),l(u),\beta_{i(u)})\nu_{i(u)}(T)(du,d\beta_{i(u)})+\\
&\displaystyle \int_{t}^{T}\int_{\mathcal{O}}h_0(u,l(u),\vartheta)\nu_0(T)(du,d\vartheta)+g_{i(T)}\big(x(T),\ell(T)\big)~~\Big]\\
&\leq u_i(t,x,l).
\end{eqnarray*}
In conclusion, the last inequality shows that $u$ is upper semi continuous, and that completes the proof.
\end{proof}
We are able now to prove one of the main results of this work, Theorem \ref{th unique caracterization for u}, dealing with;\\ \textbf{The unique characterization of the value function of the stochastic scattering problem, with optimal diffraction probability measure at the vertex selected from the own local time.}\\
\textbf{Proof of Theorem \ref{th unique caracterization for u}}:\begin{proof} We focus on showing that the value function $u$ is a continuous (Proposition \ref{pr: continuite fonction valeur}) viscosity solution of the HJB system \eqref{eq PDE Walsh} at the junction point $\bf 0$. Indeed, outside the junction point, one can use classical arguments, after a localization that pushes the spider process to be strictly included in the ray is starting from.\\
Let $h>0$ and fix $\Q^{t,x,i,l}_{\beta,\mathcal{O}}\in \mathcal{A}\big(t,x,i,l\big)$. We introduce the following stopping time:
\begin{eqnarray}\label{temp arret}
\tau_h=\inf\big\{~s\ge t_\star,~~x(s)\ge h~\big\},~~\Q^{t,x,i,l}_{\beta,\mathcal{O}}~~\text{a.s.} \end{eqnarray}
This important fact has been proven in \cite{Spider 2} Proposition 8.1, and will be used in the sequel:
\begin{eqnarray}\label{estimee centrale}
Ch^2\ge \E^{\Q^{t,x,i,l}_{\beta,\mathcal{O}}}[\tau_h-t_\star]\ge 0,~~\lim_{h \searrow 0}\frac{1}{h}\E^{\Q^{t,x,i,l}_{\beta,\mathcal{O}}}[l(\tau_h)-l_\star]=1. \end{eqnarray}
where $C>0$ is a uniform constant (depending on the data of the system) independent of $h$.\\
\textbf{Sub solution at $\bf 0$:} Let $\phi \in \mathcal{C}^{1,2,0}_{{\bf 0},1}\big(\mathcal{D}_T\big)$ be a test function  such that the local maximum of $u-\phi$ is reached at some point:
$(t_\star,0,l_\star)\in (0,T)\times [I] \times (0,+\infty)$, satisfying $(u-\phi)(t_\star, 0,l_\star)=0$. Now we use the dynamic programming principle Theorem \ref{th: PPD}, the definition of the value function $u$ given in Definition \ref{def fonction valeur} - \eqref{eq : value function}, and the fact that $u$ attains its supremum (Theorem \ref{th: compact}) for some $\Q^{t,x,i,l}_{\beta,\mathcal{O}}\in \mathcal{A}\big(t,x,i,l\big)$, to obtain:
\begin{align*}
 &\nonumber u(t_\star, 0,l_\star)=\phi(t_\star, 0,l_\star)=\sup \Big\{~~\mathbb{E}^{\P^{t,x,i,l}_{\beta,\vartheta}}\Big[~~u_{i(\tau_h)}\big(\tau,x(\tau_h),l(\tau_h)\big)~~+\\
\nonumber&\hspace{4 cm}\displaystyle\int_{t}^{\tau_h}\int_{\mathcal{B}_{i(u)}}h_{i(u)}(u,x(u),l(u),\beta_{i(u)})\nu_{i(u)}(T)(du,d\beta_{i(u)})+\\
&\nonumber \hspace{4 cm} \displaystyle \int_{t}^{\tau_h}\int_{\mathcal{O}}h_0(u,l(u),\vartheta)\nu_0(T)(du,d\vartheta)~~\Big],\\
&\hspace{4 cm} \P^{t,x,i,l}_{\beta,\vartheta}\in\mathcal{A}\big(t,x,i,l\big)~~\Big\}~\leq\\
&\mathbb{E}^{\Q^{t,x,i,l}_{\beta,\mathcal{O}}}\Big[~~\phi_{i(\tau_h)}\big(\tau,x(\tau_h),l(\tau_h)\big)+\displaystyle\int_{t}^{\tau_h}\int_{\mathcal{B}_{i(u)}}h_{i(u)}(u,x(u),l(u),\beta_{i(u)})\nu_{i(u)}(T)(du,d\beta_{i(u)})\\
&\nonumber \hspace{4 cm} +\displaystyle \int_{t}^{\tau_h}\int_{\mathcal{O}}h_0(u,l(u),\vartheta)\nu_0(T)(du,d\vartheta)~~\Big].
\end{align*}
Is follows then from the martingale property that:
\begin{align*}
&0~\leq~\mathbb{E}^{\Q^{t,x,i,l}_{\beta,\mathcal{O}}}\Big[~~\displaystyle\int_{t}^{\tau_h}\int_{\mathcal{B}_{i(u)}}\big[\partial_t\phi_{i(u)}(u,x(u),l(u))\\
&+\frac{1}{2}\sigma^2_{i(u)}(u,x(u),l(u),\beta_{i(u)})\partial_x^2\phi_{i(u)}(u,x(u),l(u))+b_{i(u)}(u,x(u),l(u),\beta_{i(u)})\times\\
&\partial_x\phi_{i(u)}(u,x(u),l(u))+h_{i(u)}(u,x(u),l(u),\beta_{i(u)})\big]\nu_{i(u)}(T)(du,d\beta_{i(u)})~~\Big]\\
&+~\displaystyle \mathbb{E}^{\Q^{t,x,i,l}_{\beta,\mathcal{O}}}\Big[~~\int_{t}^{\tau_h}\int_{\mathcal{O}}\big[\partial_l\phi(u,0,l(u))+\sum_{i=1}^I\mathcal{S}_i(u,l(u),\vartheta)+\\
&h_0(u,l(u),\vartheta)\big]\nu_0(T)(du,d\vartheta)~~\Big].
\end{align*}
We are going now to divide the last equation by $h$ and send $h \searrow 0$. For the first term, remark that assumptions $(\mathcal{H})$ together with \eqref{estimee centrale} imply:
\begin{align*}
&\overline{\lim}_{h\searrow 0}\frac{1}{h}\Big|~\mathbb{E}^{\Q^{t,x,i,l}_{\beta,\mathcal{O}}}\Big[~~\displaystyle\int_{t}^{\tau_h}\int_{\mathcal{B}_{i(u)}}\big[\partial_t\phi_{i(u)}(u,x(u),l(u))+\\
&\frac{1}{2}\sigma^2_{i(u)}(u,x(u),l(u),\beta_{i(u)})\partial_x^2\phi_{i(u)}(u,x(u),l(u))+b_{i(u)}(u,x(u),l(u),\beta_{i(u)})\times\\ &\partial_x\phi_{i(u)}(u,x(u),l(u))+h_{i(u)}(u,x(u),l(u),\beta_{i(u)})\big]\nu_{i(u)}(T)(du,d\beta_{i(u)})~~\Big]~\Big|\\
&\leq~\underline{\lim}_{h\searrow 0}~Ch~=~0, 
\end{align*}
where $C>0$ is constant independent of $h$. For the second term, write:
\begin{eqnarray}
&\nonumber \ds \frac{1}{h}\mathbb{E}^{\Q^{t,x,i,l}_{\beta,\mathcal{O}}}\Big[~~\int_{t}^{\tau_h}\int_{\mathcal{O}}\big[\partial_l\phi(u,0,l(u))+\sum_{i=1}^I\mathcal{S}_i(u,l(u),\vartheta)+\\
&\nonumber h_0(u,l(u),\vartheta)\big]\nu_0(T)(du,d\vartheta)~~\Big]\leq \\
&\nonumber\Big[\partial_l\phi(t_\star,0,l_\star)+
\underset{\vartheta\in \mathcal{O}}{\sup}\Big\{~\sum_{i=1}^I\mathcal{S}_i(t_\star,l_\star,\vartheta)\\
&\label{eq sur sol}+h_0(t_\star,l_\star,\vartheta)\Big\}\Big]\times \frac{1}{h}\mathbb{E}^{\Q^{t,x,i,l}_{\beta,\mathcal{O}}}\big[l(\tau_h)-l_\star\big]~+~\varepsilon(h).
\end{eqnarray}
The error term $\varepsilon(h)$ is given by:
\begin{eqnarray}\label{eq term erreur}
\varepsilon(h)=\frac{C}{h}\mathbb{E}^{\Q^{t,x,i,l}_{\beta,\mathcal{O}}}\big[(l(\tau_h)-l_\star)^2+(l(\tau_h)-l_\star)(\tau_h-t_\star)\big],
\end{eqnarray}
with the aid of assumption $(\mathcal{H})$ ($C>0$ stands for a uniform constant independent of $h$) and tends to $0$ as soon as $h\searrow 0$ using \eqref{estimee centrale} and the modulus of continuity of the local time given in Proposition \ref{pr: ineq}.
Using once again \eqref{estimee centrale}, it follows from \eqref{eq sur sol} that we will obtain as soon as $h\searrow 0$:
$$\partial_l\phi(t_\star,0,l_\star)+
\underset{\vartheta\in \mathcal{O}}{\sup}\Big\{~\sum_{i=1}^I\mathcal{S}_i(t_\star,l_\star,\vartheta)\partial_x\phi_i(t_\star,0,l_\star)+h_0(t_\star,l_\star,\vartheta)\Big\}~\ge~0,$$
and that completes the first \textbf{Step 1} for the sub solution case.\\
\textbf{Super solution at $\bf 0$:} Let $\phi \in \mathcal{C}^{1,2,0}_{{\bf 0},1}\big(\mathcal{D}_T\big)$ be a test function  such that the local minimum of $u-\phi$ is reached at some point:
$(t_\star,0,l_\star)\in (0,T)\times [I] \times (0,+\infty)$, satisfying $(u-\phi)(t_\star, 0,l_\star)=0$.
We argue by contradiction assuming that there exists $\varepsilon>0$ such that:
\begin{eqnarray}\label{eq contra sur sol}
\partial_l\phi(t_\star,0,l_\star)+
\underset{\vartheta\in \mathcal{O}}{\sup}\Big\{~\sum_{i=1}^I\mathcal{S}_i(t_\star,l_\star,\vartheta)\partial_x\phi_i(t_\star,0,l_\star)+h_0(t_\star,l_\star,\vartheta)\Big\}~>\varepsilon.    
\end{eqnarray}
From assumption $(\mathcal{H})$, we know that there exists:
$$\vartheta^*(t,l)=\vartheta^*\Big(t,l,\big(\partial_x\phi_i(t,0,l)_{i\in[I]}\big),\big(\mathcal{S}_{i\in[I]}\big),h_0\Big),$$
such that $\forall (t,l)\in[0,T]\times[0,+\infty)$:
\begin{align*}
 \underset{\vartheta\in \mathcal{O}}{\sup}\Big\{~\sum_{i=1}^I\mathcal{S}_i(t,l,\vartheta)\partial_x\phi_i(t,0,l)+h_0(t,l,\vartheta)\Big\}=\\\sum_{i=1}^I\mathcal{S}_i(t,l,\vartheta^*(t,l))\partial_x\phi_i(t,0,l)+h_0(t,l,\vartheta^*(t,l)).
\end{align*}
Once again from assumption $(\mathcal{H})$, we have that the coefficients of diffraction are uniformly Lipschitz continuous w.r.t to the control variables $\vartheta \in \mathcal{O}$. We can use the central Theorem 2.2 of \cite{Spider}, and build some $\Q^{t,x,i,l}_{\beta,\mathcal{O}}\in \mathcal{A}\big(t,x,i,l\big)$ possessing exactly the spinning measure $(t,l)\to \mathcal{S}_i(t,l,\vartheta^*(t,l)).$ Now we proceed like in the sub solution case to obtain:
(using the dynamic programming principle Theorem \ref{th: PPD}, and the definition of the value function $u$ given in Definition \ref{def fonction valeur} - \eqref{eq : value function})  
\begin{align*}
 &\nonumber u(t_\star, 0,l_\star)=\phi(t_\star, 0,l_\star)=\sup \Big\{~~\mathbb{E}^{\P^{t,x,i,l}_{\beta,\vartheta}}\Big[~~u_{i(\tau_h)}\big(\tau,x(\tau_h),l(\tau_h)\big)~~+\\
\nonumber&\hspace{4 cm}\displaystyle\int_{t}^{\tau_h}\int_{\mathcal{B}_{i(u)}}h_{i(u)}(u,x(u),l(u),\beta_{i(u)})\nu_{i(u)}(T)(du,d\beta_{i(u)})+\\
&\nonumber \hspace{4 cm} \displaystyle \int_{t}^{\tau_h}\int_{\mathcal{O}}h_0(u,l(u),\vartheta)\nu_0(T)(du,d\vartheta)~~\Big],\\
&\hspace{4 cm} \P^{t,x,i,l}_{\beta,\vartheta}\in\mathcal{A}\big(t,x,i,l\big)~~\Big\}~\ge\\
&\mathbb{E}^{\Q^{t,x,i,l}_{\beta,\mathcal{O}}}\Big[~~\phi_{i(\tau_h)}\big(\tau,x(\tau_h),l(\tau_h)\big)+\displaystyle\int_{t}^{\tau_h}\int_{\mathcal{B}_{i(u)}}h_{i(u)}(u,x(u),l(u),\beta_{i(u)})\nu_{i(u)}(T)(du,d\beta_{i(u)})\\
&\nonumber \hspace{4 cm} +\displaystyle \int_{t}^{\tau_h}\int_{\mathcal{O}}h_0(u,l(u),\vartheta)\nu_0(T)(du,d\vartheta)~~\Big].
\end{align*}
Dividing by $h>0$ and using the martingale property, and once the central estimates \eqref{estimee centrale}, we claim that:
\begin{align}
&0\ge \nonumber \frac{1}{h}\mathbb{E}^{\Q^{t,x,i,l}_{\beta,\mathcal{O}}}\Big[~~\int_{t}^{\tau_h}\int_{\mathcal{O}}\big[\partial_l\phi(u,0,l(u))+\sum_{i=1}^I\mathcal{S}_i(u,l(u),\vartheta^*(u,l(u)))+\\
&\nonumber h_0(u,l(u),\vartheta^*(u,l(u)))\big]\nu_0(T)(du,d\vartheta)~~\Big]-Ch\ge \\
&\nonumber\Big[\partial_l\phi(t_\star,0,l_\star)+
\underset{\vartheta\in \mathcal{O}}{\sup}\Big\{~\sum_{i=1}^I\mathcal{S}_i(t_\star,l_\star,\vartheta)\\
&+h_0(t_\star,l_\star,\vartheta)\Big\}\Big]\times \frac{1}{h}\mathbb{E}^{\Q^{t,x,i,l}_{\beta,\mathcal{O}}}\big[l(\tau_h)-l_\star\big]~-~\varepsilon(h)-Ch.
\end{align}
where $C>0$ is constant independent of $h$, and the error term $\varepsilon(h)$ is given like in \eqref{eq term erreur} and tends to $0$ as soon as $h\searrow 0$. Therefore, sending $h\searrow 0$, from \eqref{estimee centrale} and \eqref{eq contra sur sol}, we obtain:
$$0\ge \varepsilon >0,$$
and that leads to a contradiction and completes the proof.
\end{proof}
\section{On stochastic scattering control problems without dependency w.r.t to the local-time variable.}\label{sec control without l}
We add here a short section, in which we will discuss briefly on stochastic scattering control problem without dependency w.r.t to the local-time variable. This type of problem appears as a particular case of the previous results studied in this work, when the optimal diffraction probability measure is selected from the own local time at the vertex $\bf 0$.

Let us briefly explain the main lines that one has to undertake to solve such problems.\\
\textbf{ In the sequel,  we assume the assumption $(\mathcal{H})$ and that all the data $(\mathcal{D})$ have no dependence with respect to the variable $l$.}\\
-a) First, one has to introduce and study the uniqueness of the following parabolic HJB system posed on the star-shaped network $\mathcal{N}$:
\begin{eqnarray}\label{eq PDE Walsh sans l}
&\begin{cases}
\textbf{HJB equation on each ray}:\\
\partial_t u_i(t,x)+\underset{\beta_i\in \mathcal{B}_i}{\sup}\Big\{\ds \frac{1}{2}\sigma_i^2(t,x,\beta_i)\partial^2_xu_i(t,x)+\\
b_i(t,x,\beta_i)\partial_xu_i(t,x)+h_i(t,x,\beta_i)\Big\}=0,~~(t,x)\in(0,T)\times \R_+,\\
\textbf{Non linear Kirchhoff's boundary transmission at } \bf 0: \\
\underset{ \vartheta \in \mathcal{O}}{\sup} \Big\{\ds\sum_{i=1}^I\mathcal{S}_i(t,\vartheta)\partial_xu_i(t,0)+h_0(t,\vartheta)\Big\}=0,~~t\in(0,T),~~\\
\textbf{Terminal condition}:\\
u_i(T,x)=g_i(x),~~\forall x\in\R_+.~~\forall i\in[I],\\
\textbf{Continuity condition at } \bf 0:\\
\forall (i,j)\in[I]^2,~~\forall t\in \R_+,~~u_i(t,0)=u_j(t,0).
\end{cases}
\end{eqnarray}
We claim that the comparison theorem holds true for continuous viscosity solutions of the last system \eqref{eq PDE Walsh sans l}. To be convinced, one can follow exactly the same lines of proof of Theorem 2.4 given in \cite{Ohavi Walsh PDE}, using a small exponential perturbation involving the deterministic local-time variable $l$ and the comparison Theorem \ref{th compa visco} of this work.\\
-b) Secondly, the construction of optimal rules $\P^{t,x,i}_{\beta,\mathcal{O}}$ can be achieved with the same arguments given in the work \cite{Spider}. Note that in this case, the canonical space should not be extended and we can consider only the space of continuous functions living in $\mathcal{N}$, namely:  $\Phi=\mathcal{C}^{\mathcal{N}}\big([0,T]\big)$. The construction can be done then using an adapted approximate sequence built by concatenation of martingale problems, and a tension argument. For the uniqueness, we can also proceed like in the proof of Theorem 3.1.II of \cite{Spider}, but without the use of the system involving the local-time Kirchhoff's boundary condition, for which the well-posedness has been obtained in \cite{Martinez-Ohavi EDP}. It is enough to use the well-posedness of linear parabolic system posed on star-shaped networks, studied in Section 3 of \cite{Martinez-Ohavi EDP}. Therein, a classical Kirchhoff's boundary condition is considered, and the results also improve the original ones of \cite{Von-Below}. Notably, Theorem 3.5 of \cite{Martinez-Ohavi EDP} states the unique solvability with weaker assumptions, for the compatibility conditions and the regularity of the coefficients of diffraction, that are more natural when we link this type of PDE's system to the generators of spider diffusions.\\
-c) Finally, it follows that in this case when there is no dependency with the local-time variable in the spinning measure and coefficients of diffusion, the same ideas will lead to obtain the  weak dynamic programming principle Theorem \ref{th: PPD}. On the other side, the unique characterization of the value function $u$ given in \eqref{eq : value function} (without dependency w.r.t to the local-time) can be achieved with the aid of system \eqref{eq PDE Walsh sans l}, and the same arguments of the proof of Theorem \ref{th unique caracterization for u}.
\appendix 
\section{Some analysis tools}\label{app : element analyse funct}
We recall here some definitions and functional analysis tools that were used in this contribution, essentially to obtain the compactness criterion Theorem \ref{th : compactV}, for the admissible rules $\mathcal{V}\big([0,T]\times \mathcal{O}\big)$ at the vertex $\bf 0$. Most of the following results are contained in \cite{jj}. 

In the sequel, we consider:
\begin{eqnarray*}
\begin{cases}
(X,T)\text{ a topological space and } \Sigma \text{ a } \sigma \text{ algebra on } X,\\
(E,\mathcal{E}) \text{ a measurable space,}\\
(F,d) \text{ be a Polish space, endowed with its metric } d  \text{ and } \mathbb{B}(F)\text{ its Borel  algebra.}  
\end{cases}
\end{eqnarray*}
\begin{Definition}
$(E,\mathcal{E})$ is said to be countably generated, if there exists a countable base generating $\mathcal{E}$. Namely there exists a sequence $O_{n}$ of $\mathcal{E}$, such that $\mathcal{E}=\sigma(O_n,n\in \mathbb{N})$.
\end{Definition}
Recall that since $F$ is a Polish space, the measurable space $(F,\mathbb{B}(F))$ is countably generated, (see for instance Proposition 3.1 in \cite{Ch P}).
\begin{Definition}
Let $P$ be a measure on $(X, \Sigma)$. We say that $P$ is regular if for any measurable subset $B\in\Sigma$ 
\begin{eqnarray*}
&P(B)=\sup\Big\{P(K),~K \text{ closed, }K\in \Sigma,~K\subset B  \Big\}\\
&=\inf\Big\{P(O),~O \text{ open, } O\in \Sigma,~B\subset O \Big\}.
\end{eqnarray*}
\end{Definition}
Recall that any Borel  probability measure, or in other terms any probability measure on a metric space endowed with its $\sigma$-Borel  algebra, is regular (see for instance Proposition 2.3 in \cite{Ch P}).\\
In the rest of this Appendix, we denote by:\\
-$L^{\infty}(E)$ the set consisting of all measurable real bounded maps on $(E,\mathcal{E})$.\\
-$\mathcal{C}(F)$, (resp. $\mathcal{C}_u(F)$), are the set of continuous (resp. uniformly continuous) bounded functions on $F$.\\
-$L^{\infty}(E\times F)$ is the set of measurable bounded real functions defined on $\Big(E\times F,\mathcal{E}\otimes \mathbb{B}(F)\Big)$.\\
-$\mathcal{M}(E)$ the set consisting of nonnegative finite measures on $(E,\mathcal{E})$.\\
-$\mathcal{M}(F)$ the set consisting of nonnegative finite measures on $(F,\mathbb{
B}(F))$.\\
-$\mathcal{M}(E\times F)$
the set consisting of nonnegative finite measures on $\Big(E\times F,\mathcal{E}\otimes \mathbb{B}(F)\Big)$.\\
We introduce furthermore: 
\begin{eqnarray*}
&L^{\infty}_{mc}(E\times F):=\Big\{f\in \mathcal{L}^\infty(E\times F),~\text{s.t.}~x\mapsto f(s,x)\in \mathcal{C}(F),~~\forall s\in E\Big\},\\
&L^{\infty,1}_{mc}(E\times F):=\Big\{f\in L^\infty_{mc}(E\times F),~\text{s.t.}~\exists A\in \mathcal{E},~~g\in \mathcal{C}_u(F),~~f(x,z)=\mathbf{1}_{A}(x)g(z)\Big\},\\
&L^{\infty,2}_{mc}(E\times F)):=\Big\{f\in L^\infty_{mc}(E\times F),~\text{s.t.}~\exists (A_n) \text{ a partition of } \mathcal{E},\\ 
&\text{ and a sequence } (g_n) \text{ of } \in\mathcal{C}_u(F),~\text{s.t.}~f(x,z)=\displaystyle \sum_n \mathbf{1}_{A_n}(x)g_n(z).\Big\}.
\end{eqnarray*}
We will also denote
$\mathcal{M}(E)$ (resp. $\mathcal{M}(F)$, $\mathcal{M}(E\times F)$) by $\mathcal{M}_m(E)$, (resp.$\mathcal{M}_c(F)$, $\mathcal{M}_{mc}(E\times F)$) when they are endowed with the finest topology making continuous the following family of linear forms $(\theta_f)_{f\in L^{\infty}(E)}$, defined by 
\begin{eqnarray*}
\theta_f:\begin{cases}
\mathcal{M}_m(E)\to \R \\
\nu\mapsto \nu(f)=\displaystyle \int_E f d\nu
\end{cases}.
\end{eqnarray*}
(resp. $(\theta_f)_{f\in \mathcal{C}(F)}$
\begin{eqnarray*}
\theta_f:\begin{cases}
\mathcal{M}_c(F)\to \R \\
\nu\mapsto \nu(f)=\displaystyle \int_F f d\nu
\end{cases},
\end{eqnarray*}
$(\theta_f)_{f\in L^{\infty}_{mc}(E\times F)}$,
\begin{eqnarray*}
\theta_f:\begin{cases}
\mathcal{M}_{mc}(E\times F)\to \R \\
\nu\mapsto \nu(f)=\displaystyle \int_{E\times F} f d\nu
\end{cases}.)
\end{eqnarray*}
We identify $\mathcal{M}_{mc}(E\times F)$ (resp. $\mathcal{M}_{m}(E)$, $\mathcal{M}_{c}(F)$), as subsets of the dual spaces $L^\infty_{mc}(E\times F)^{'}$ (resp. $L^\infty_{m}(E)^{'}$, $\mathcal{C}(F)^{'})$, endowed with the weak topologies $*\sigma\Big(L^\infty_{mc}(E\times F)^{'},L^\infty_{mc}(E\times F)\Big)$ (resp. $*\sigma\Big(L^\infty_{m}(E)^{'},L^\infty_{m}(E)\Big)$, $*\sigma\Big(\mathcal{C}(F)^{'},\mathcal{C}(F)\Big)$).\\
We recall that a sequence $\nu_n$ of $L^\infty_{mc}(E\times F)^{'}$ (resp. $L^\infty_{m}(E)^{'}$, $\mathcal{C}(F)^{'})$, converges to $\nu \in L^\infty_{mc}(E\times F)^{'}$, (resp. $L^\infty_{m}(E)^{'}$, $\mathcal{C}(F)^{'})$ for the weak topology $*$, and we denote  $\nu_{n}\overset{*}{\rightharpoonup} \nu$, if and only if
\begin{eqnarray*}
&\forall f\in L^\infty_{mc}(E\times F),~~\nu_n(f) \xrightarrow[]{n\to +\infty} \nu(f),\\&
\Big(~~\text{resp.}~~\forall f\in L^\infty_{m}(E),~~\nu_n(f) \xrightarrow[]{n\to +\infty} \nu(f),~~\forall f\in \mathcal{C}(F),~~\nu_n(f) \xrightarrow[]{n\to +\infty} \nu(f)~~\Big).
\end{eqnarray*}
For any $\nu \in \mathcal{M}(E\times F)$, we denote by $\nu^E$ (resp. $\nu^F$), the marginal of $\nu$ on $E$ (resp. on $F$), defined by
\begin{eqnarray*}
\nu^E(dx)~~~=~~\int_{z\in F}\nu(dz),~~~~\nu^F(dz)~~~=~~\int_{x\in E}\nu(dx).
\end{eqnarray*}
\begin{Proposition}\label{pr : metri topo faible}
Suppose that $E$ is countably generated, then $\mathcal{M}_{mc}(E\times F)$ is metrizable. (See Proposition 2.10 in \cite{jj}.)
\end{Proposition}
\begin{Theorem}\label{th : compa rela}
Let $\mathcal{N}$ be a subset of $\mathcal{M}_{mc}(E\times F)$. Then $\mathcal{N}$ is relatively compact if and only if \\
(i) ~$\Big\{\nu^F,~\nu \in \mathcal{N}\Big\}$ is relatively compact in $\mathcal{M}_{m}(E)$,\\
(ii) $\Big\{\nu^E,~\nu \in \mathcal{N}\Big\}$ is relatively compact in $\mathcal{M}_{c}(F)$.\\
(See Proposition 2.10 in \cite{jj}.)
\end{Theorem}
\begin{Theorem}\label{th : Riez}
Let $\phi$ be a positive linear form defined on the vector space generated by $L^{\infty,1}_{mc}(E\times F)$ such that:\\
(i) the following map $$\begin{cases}
\mathcal{E}\to \R \\
A\mapsto \phi(\mathbf{1}_A\otimes 1)
\end{cases}$$
is a measure on $(E,\mathcal{E})$, where we have set for each $(x,z)\in E\times F$, $\mathbf{1}_A\otimes 1(x,z)=1$, if $x\in A$ and $\mathbf{1}_A\otimes 1(x,z)=0$, if $x\notin A$.\\
(ii) for each $\varepsilon>0$, there exists a compact set $K_\varepsilon$ of $F$ such that $\phi(1)-\phi(1\otimes f)\leq \varepsilon$, for any $f\in \mathcal{C}_u(F)$, satisfying $\mathbf{1}_{K_\varepsilon}\leq f\leq 1$, where we have set for each $(x,z)\in E\times F$, $1(x,z)=1$, and $1\otimes f(x,z)=f(z)$.\\
If $\phi$ satisfies $(i)-(ii)$, then there exists $\nu\in \mathcal{M}_{mc}(E\times F)$ such that
\begin{eqnarray*}
\forall f\in L^{\infty,1}_{mc}(E\times F),~~\phi(f)=\int_{E\times F}fd\nu.
\end{eqnarray*}
(See Theorem 2.6 in \cite{jj}).
\end{Theorem}
\begin{Lemma}\label{lm : densité }
Let $K$ be a compact set of $F$, and $f\in L^\infty_{mc}(E\times F)$. Then there exists a sequence $f_n$ of $L^{\infty,2}_{mc}(E\times F)$ converging to $f$ uniformly on $E\times K$.
(See Lemma 2.5 in \cite{jj}).
\end{Lemma}
We end this appendix by stating a lemma of measurable selection in the spirit of \cite{ElKaroui}, and some intermediate results that are technical, but of high importance if one want to adapt the ideas of \cite{ElKaroui} to prove the dynamic programming principle, adapted to our purpose Theorem \ref{th: PPD}.
\begin{Corollary}\label{cr : measurabilité graph fermé}(See Corollary 5.2.2 in \cite{ElKaroui}). Let $Y$ be a separable metric space. Let $y\mapsto K_y$ be a map from $Y$ into $\text{comp}(X)$, the set of the compact subsets of a separable metric space $X$, (where $\text{comp}(X)$ is endowed with the Hausdorff topolgy). If the map $y\mapsto K_y$
is closed: for any sequence $y_n$ converging to $y$, $x_n=K_{y_n}$ has a limit point in $K_y$, then the map $y\mapsto K_y$ is Borel.
\end{Corollary}
\begin{Proposition}\label{cr: sele mesur} (see for instance Corollary 5.4 in \cite{ElKaroui}). 
Let $\mathcal{G},\mathcal{H}$ be two separable metric spaces. Let $\phi$ a lower semi continuous real function on $\mathcal{G}\times \mathcal{H}$ and $h\mapsto K_{h}$ a measurable map from $\mathcal{H}$ into $\text{comp}(\mathcal{G})$, (the set of compacts sets of $\mathcal{G}$, endowed with the Hausdorff metric). Then  \\
-the map : $v(h):=\inf\Big\{\phi(g,h),g \in K_{h}\Big\}$ is a Borel  function and $h\mapsto M_{h} := \Big\{g,v(h)=\phi(h,g),g\in K_{h}\Big\}$ is a measurable map of $\mathcal{H}$ into $\text{comp}(\mathcal{G})$.
\\
-for each probability measure $P$ on $\mathcal{H}$: 
\begin{eqnarray*}
&\displaystyle\int v(h)dP(h)~~=~~\int\inf\Big\{~\phi(g,h),~g \in K_{h}~\Big\}dP(h)\\
&=~~\inf\Big\{~~\displaystyle\int \phi(\beta(h),h)dP(h),~\beta:\mathcal{H}\to \mathcal{G},~\text{measurable},~\beta(h) \in K_{h}~\Big\}.
\end{eqnarray*}
\end{Proposition}
\begin{Proposition}\label{pr : countably gener sub algebra}
Let $\tau$ be a $(\Phi)_{t\leq s\leq T}$ stopping time, then: 
$$\Phi_\tau=\sigma\Big(X(s\wedge\tau),s\leq T\Big),$$
and $\Phi_\tau$ is countably generated. 
\end{Proposition}
\begin{proof}
Recall that
\begin{eqnarray*}
\Phi_{\tau} = \Big\{B\in \Phi_{T},~ B\cap\{\tau\leq t\}\in \Phi_{t},~\forall t\in[0,T]\Big\}.
\end{eqnarray*}
Recall that the space where is defined our canonical process $X(\cdot)$, is giving by:
\begin{eqnarray*}
\Phi~~=~~\mathcal{C}^{\mathcal{N}}([0,T])\times \Big(\prod_{i=1}^I \mathcal{U}([0,T]\times \mathcal{B}_i)\Big)\times \mathcal{V}([0,T]\times \mathcal{O}),
\end{eqnarray*}
is Polish.\\
We can use then the same arguments of the proof of Lemma 1.3.3 in \cite{Stroock}, to conclude.
\end{proof}

\end{document}